\documentclass[aop,citesort,MSNbibl,seceqn,dvips]{arximspdf}

\doi{10.1214/10-AOP611}
\volume{40}
\issue{1}
\pubyear{2012}
\firstpage{213}
\lastpage{244}

\makeatletter

\newcommand{\fraca}[2]{{#1/#2}}
\newcommand{\fracb}[2]{{(#1)/#2}}
\def\sF{{\mathcal F}}
\def\bD{{\mathbb D}}
\def\bP{{\mathbb P}}
\def\bR{{\mathbb R}}
\def\wt{\widetilde}
\def\wh{\widehat}
\def\E{{\mathbb E}}
\def\P{{\mathbb P}}
\def\eps{\varepsilon}
\def\wh{\widehat}
\newcommand{\cal}{\mathcal}
\newcommand{\eqref}[1]{(\ref{#1})}
\newtheorem{theo}{Theorem}[section]
\newproclaim{remark}[theo]{Remark}
\newtheorem{lemma}[theo]{Lemma}
\newtheorem{prop}[theo]{Proposition}
\newtheorem{cor}[theo]{Corollary}
\makeatother

\begin{document}
\begin{frontmatter}

\title{Sharp heat kernel estimates for relativistic stable processes in open sets}
\runtitle{Heat kernel estimates for relativistic stable processes}

\begin{aug}
\author[A]{\fnms{Zhen-Qing} \snm{Chen}\thanksref{t1}\ead[label=e1]{zchen@math.washington.edu}},
\author[B]{\fnms{Panki} \snm{Kim}\corref{}\thanksref{t2}\ead[label=e2]{pkim@snu.ac.kr}}
\and
\author[C]{\fnms{Renming} \snm{Song}\ead[label=e3]{rsong@math.uiuc.edu}}
\runauthor{Z.-Q. Chen, P. Kim and R. Song}
\affiliation{University of Washington, Seoul National University
and~University~of~Illinois}
\address[A]{Z.-Q. Chen\\
Department of Mathematics\\
University of Washington\\
Seattle,
Washington 98195\\
USA\\
\printead{e1}} 
\address[B]{P. Kim\\
Department of Mathematical Sciences\\
Seoul National University\\
San 56-1 Shinrim-dong Kwanak-gu\\
Seoul 151-747\\
Korea\\
\printead{e2}}
\address[C]{R. Song\\
Department of Mathematics\\
University of Illinois\\
Urbana, Illinois 61801\\
USA\\
\printead{e3}}
\end{aug}
\thankstext{t1}{Supported in part
by NSF Grants DMS-06-00206 and DMS-09-06743.}
\thankstext{t2}{Supported by Basic Science Research Program through
the National
Research Foundation of Korea (NRF) grant funded by the Korean
government (MEST) (2010-0028007).}

\received{\smonth{2} \syear{2010}}
\revised{\smonth{8} \syear{2010}}

%
\begin{abstract}
In this paper, we establish sharp two-sided estimates for the
transition densities of relativistic stable processes [i.e., for the
heat kernels of the operators $m-(m^{2/\alpha}-\Delta)^{\alpha/2}$]
in $C^{1, 1}$ open sets.
Here $m>0$ and $\alpha\in(0, 2)$.
The estimates are uniform in $m\in(0, M]$
for each fixed $M>0$. Letting $m\downarrow0$, we recover the
Dirichlet heat kernel estimates for
$\Delta^{\alpha/2}:= - (-\Delta)^{\alpha/2}$
in
$C^{1,1}$ open sets obtained in~\cite{CKS}. Sharp two-sided
estimates are also obtained for Green functions of relativistic
stable processes in bounded $C^{1,1}$ open sets.
\end{abstract}

%

\end{frontmatter}

\section{Introduction}\label{intro}

Throughout this paper we assume that $d\ge1$ and $\alpha\in(0,
2)$. For any $m>0$,
a relativistic $\alpha$-stable process $X^m$
on $\bR^d$ with mass $m$ is a L\'{e}vy process with characteristic
function given by
%
%
\begin{equation}\label{e:ch}
\E\bigl[\exp\bigl(i \xi\cdot(X^m_t-X^m_0)
\bigr) \bigr] = \exp\bigl(-t \bigl( (|\xi|^2+ m^{2/\alpha}
)^{\alpha/2}-m \bigr) \bigr), \qquad\xi\in\bR^d.\hspace*{-35pt}
\end{equation}
The limiting case $X^0$, corresponding to $m=0$, is a (rotationally)
symmetric $\alpha$-stable (L\'{e}vy) process on $\bR^d$ which we will
simply denote as $X$. The infinitesimal generator of $X^m$ is
$m-(m^{2/\alpha} -\Delta)^{\alpha/2}$. Note that when $m=1$, this
infinitesimal generator reduces to $1-(1 -\Delta)^{\alpha/2}$. Thus
the $1$-resolvent kernel of the relativistic $\alpha$-stable
process $X^1$ on $\bR^d$ is just the Bessel potential kernel. (See
\cite{BMR} for more on this connection.) When $\alpha=1$, the
infinitesimal generator reduces to the so-called free relativistic
Hamiltonian $ m - \sqrt{-\Delta+ m^{2}}$. The operator $m -
\sqrt{-\Delta+ m^{2}}$ is very important in mathematical physics
due to its application to relativistic quantum mechanics. Physical
models related to this operator have been much studied over the past
30 years, and there exists a huge literature on the properties of
relativistic Hamiltonians (see, e.g.,
\cite{CMS,FL,He,Lieb,LY,W}). For recent papers in the mathematical
physics literature related to the relativistic Hamiltonian, we refer
the readers to~\cite{DOS,FLS1,FLS2,S} and the references therein.
Various fine properties of relativistic $\alpha$-stable processes
have been studied recently in~\cite{BMR,CK2,CS4,GR2,K2,KL,KS,R}.

The objective of this paper is to establish (quantitatively) sharp
two-sided estimates on the transition density $p^m_D(t, x, y)$ of
the subprocess of $X^m$ killed upon exiting any $C^{1,1}$ open set
$D\subset\bR^d$. The density function $p^m_D(t, x, y)$ is also the
heat kernel of the restriction of $m-(m^{2/\alpha} -\Delta)^{\alpha
/2}$ in $D$ with zero exterior condition. Recall that an open set
$D$ in $\bR^d$ (when $d\ge2$) is said to be a (uniform) $C^{1,1}$
open set if there are (localization radius) $ R>0 $ and
$\Lambda_0>0$ such that for every $z\in\partial D$, there exist a
$C^{1,1}$-function $\varphi=\varphi_z\dvtx \bR^{d-1}\to\bR$ satisfying
$\varphi(0)=0, \nabla\varphi(0)=(0, \ldots, 0)$, $| \nabla\varphi
(x)-\nabla\varphi(z)| \leq\Lambda_0 |x-z|$ and an orthonormal
coordinate system $CS_z\dvtx y=(y_1, \ldots, y_{d-1}, y_d):=(\wt y,
y_d)$ with origin at $z$ such that $ B(z, R )\cap D= \{y=(\wt y,
y_d)\in B(0, R) \mbox{ in } CS_z\dvtx y_d
> \varphi(\wt y) \}$. By a $C^{1,1}$ open set in $\bR$ we mean an open
set which can be expressed as the union of disjoint intervals so
that the minimum of the lengths of all these intervals is positive,
and the minimum of the distances between these intervals is
positive. For $x\in\bR^d$, let $\delta_D(x)$ denote the Euclidean
distance between $x$ and $D^c$ and $\delta_{\partial D}(x)$ the
Euclidean distance between $x$ and $\partial D$. It is well known
that a $C^{1, 1}$ open set $D$ satisfies both the \textit{uniform
interior ball condition} and the \textit{uniform exterior ball
condition}: there exists $r_0 < R $ such that for every $x\in D$
with $\delta_{\partial D}(x)< r_0$ and $y\in\bR^d \setminus
\overline D$ with $\delta_{\partial D}(y)<r_0$, there are $z_x,
z_y\in\partial D$ so that $|x-z_x|=\delta_{\partial D}(x)$,
$|y-z_y|=\delta_{\partial D}(y)$ and that $B(x_0, r_0)\subset D$ and
$B(y_0, r_0)\subset\bR^d \setminus\overline D$, where
$x_0=z_x+r_0(x-z_x)/|x-z_x|$ and $y_0=z_y+r_0(y-z_y)/|y-z_y|$. In
fact, $D$ is $C^{1, 1}$ if and only if $D$ satisfies the uniform
interior ball condition and the uniform exterior ball condition
(see~\cite{AKSZ}, Lemma 2.2). In this paper we call the pair $(r_0,
\Lambda_0)$ the characteristics of the $C^{1,1}$ open set $D$.

The main result of this paper is Theorem~\ref{t:main} below. The
open set $D$ below is not necessarily bounded or connected. In this
paper, we use ``$:=$'' as a way of definition. For $a, b\in\bR$,
$a\wedge b:=\min\{a, b\}$ and $a\vee b:=\max\{a, b\}$.

\begin{theo}\label{t:main}
Suppose that $D$ is a $C^{1,1}$ open set in $\bR^d$ with
$C^{1,1}$ characteristics $(r_0, \Lambda_0)$.
\begin{longlist}[(ii)]
\item[(i)]
For any $M>0$ and $T>0$,
there exists $C_1= C_1(\alpha, r_0, \Lambda_0, M, T)>1$
such that for any $m\in(0, M]$ and
$(t,x,y) \in(0, T]\times D\times D$,
%
%
\begin{eqnarray}\label{e:1.2}
&& \frac1{C_1} \biggl( 1\wedge
\frac{\delta_D(x)^{\alpha/2}}{\sqrt{t}} \biggr)\biggl ( 1\wedge
\frac{\delta_D(y)^{\alpha/2}}{\sqrt{t}} \biggr) \biggl(
t^{-d/\alpha}
\wedge\frac{t
\phi(m^{1/\alpha} |x-y|)}
{|x-y|^{d+\alpha}} \biggr)\nonumber
\\
&& \qquad \le p^m_D(t, x, y)\\
&& \qquad \le C_1 \biggl( 1\wedge
\frac{\delta_D(x)^{\alpha/2}}{\sqrt{t}} \biggr)\biggl ( 1\wedge
\frac{\delta_D(y)^{\alpha/2}}{\sqrt{t}} \biggr)\biggl (
t^{-d/\alpha}
\wedge\frac{t
\phi(m^{1/\alpha} |x-y|/16)}
{|x-y|^{d+\alpha}} \biggr)
,
\nonumber
\end{eqnarray}
where $\phi(r)= e^{-r}(1+r^{(d+\alpha-1)/2})$.
\item[(ii)] Suppose in addition that $D$ is bounded.
For any $M>0$ and $T>0$,
there exists $C_2= C_2(\alpha, r_0, \Lambda_0, M, T, \operatorname
{diam}(D)) >1$
such that for any $m\in(0, M]$ and
$(t,x,y) \in[T, \infty)\times D\times D$,
\begin{eqnarray*}
C_2^{-1} e^{- t \lambda^{\alpha, m, D}_1 } \delta_D
(x)^{\alpha/2} \delta_D (y)^{\alpha/2} &\leq& p^m_D(t, x, y)\\
&\leq& C_2 e^{-t \lambda^{\alpha, m, D}_1 } \delta_D
(x)^{\alpha/2} \delta_D (y)^{\alpha/2} ,
\end{eqnarray*}
where $\lambda^{\alpha, m, D}_1>0$ is the smallest eigenvalue of the
restriction of $(m^{2/\alpha}-\Delta)^{\alpha/2}-m$ in $D$ with
zero exterior condition.
\end{longlist}
\end{theo}

\begin{remark}\label{R:1.2}
(i) Note that the estimates in Theorem~\ref{t:main} are
uniform in $m\in(0, M]$. When $m\downarrow0$, $m-
(m^{2/\alpha}-\Delta)^{\alpha/2}$ converges to the fractional
Laplacian $\Delta^{\alpha/2}:= - (-\Delta)^{\alpha/2}$ in the
distributional sense, and it is easy to check that $X^m$ converges
weakly to $X$ in the Skorokhod space $\bD([0, \infty), \bR^d)$. It
follows from the uniform H\"{o}lder continuity result of~\cite{CK},
Theorem 4.14, that $p^m_D(t, x, y)$ converges pointwise to $p_D(t, x,
y)$, the transition density function of the subprocess $X^D$ of $X$
in $D$. Furthermore, when $D$ is bounded, by~\cite{CS9}, Theorem 1.1,
$\lim_{m\downarrow0} \lambda^{\alpha, m, D}_1 =
\lambda^{\alpha, D}_1$, the smallest eigenvalue of
$(-\Delta)^{\alpha/2}$ in $D$ with zero exterior condition. So
letting $m\downarrow0$ in Theorem~\ref{t:main} recovers the sharp
two-sided estimates of $p_D(t, x, y)$ for $C^{1,1}$ open set $D$,
which were obtained recently in~\cite{CKS}. We emphasize here that
the proof of Theorem~\ref{t:main} of this paper uses the main
results of~\cite{CKS}, so the above remark should not be interpreted
as that passing $\alpha\to0$ gives a new proof of the main results
of~\cite{CKS}.

(ii)
When $D$ is bounded, the functions $(x, y)\mapsto\phi(m^{1/\alpha}
|x-y|/16)$ and $(x, y)\mapsto\phi(m^{1/\alpha} |x-y|)$ on
$D\times D$ are bounded between two positive constants independent
of $m \in(0, M]$. Thus it follows from Theorem~\ref{t:main}(i)
above and
\cite{CKS}, Theorem~1.1(i),
that, for each $T>0$, the
heat kernel $p^m_D(t, x, y)$ is uniformly comparable to the heat
kernel $p_D(t, x, y)$ on $(0, T]\times D\times D$ when~$D$ is a
bounded $C^{1, 1}$ open set. However when $D$ is unbounded, these
two are not comparable.

(iii)
In fact, the upper bound estimates in both Theorem~\ref{t:main} and
Theorem~\ref{T:1.2} below hold for any open set $D$ satisfying (a
weak version of) the \textit{uniform exterior ball condition} in place
of the $C^{1,1}$ condition, while the lower bound estimates in both
Theorem~\ref{t:main} and Theorem~\ref{T:1.2} below hold for any open
set $D$ satisfying the \textit{uniform interior ball condition} in
place of the $C^{1,1}$ condition. (See the paragraph before
Lemma~\ref{ub11} for the definition of the weak version of the
uniform exterior ball condition.)\vadjust{\goodbreak}

(iv) Let $p^m(t, x, y)$ denote
the transition density function for $X^m$.
Then in view of \eqref{e:2.10} and the estimates on $p^m(t, x, y)$ to be
given below in Theorem~\ref{T:lbRd}, the estimate \eqref{e:1.2} can
be restated as
%
%
\begin{eqnarray}\label{e:1.3}
&&  \frac1{C_1} \biggl( 1\wedge
\frac{\delta_D(x)^{\alpha/2}}{\sqrt{t}} \biggr) \biggl( 1\wedge
\frac{\delta_D(y)^{\alpha/2}}{\sqrt{t}} \biggr) p^m(t, x, y)
\nonumber\hspace*{-35pt}
\\[-4pt]
\\[-12pt]
&&  \qquad \leq p^m_D(t, x, y) \le C_1 \biggl(
1\wedge
\frac{\delta_D(x)^{\alpha/2}}{\sqrt{t}} \biggr)\biggl ( 1\wedge
\frac{\delta_D(y)^{\alpha/2}}{\sqrt{t}} \biggr) p^m(t, x/16, y/16).
\nonumber\hspace*{-35pt}
\end{eqnarray}
\end{remark}

Though the heat kernel estimates for symmetric diffusions (such as
Aronson's estimates) have a long history, the study of sharp
two-sided estimates on the transition densities of jump processes in
$\bR^d$ started quite recently. See~\cite{C,CKK,CK,CK2} and the
references therein. Due to the complication near the boundary, the
investigation of sharp two-sided estimates on the transition
densities of jump processes in open sets is even more recent. In
\cite{CKS}, we obtained sharp two-sided estimates for the
transition density of the symmetric $\alpha$-stable process killed
upon exiting any $C^{1, 1}$ open set $D\subset\bR^d$. That was the
first time sharp two-sided estimates were established for Dirichlet
heat kernels of nonlocal operators. Subsequently, we obtained in
\cite{CKS1} sharp two-sided heat kernel estimates for censored
stable processes in $C^{1,1}$ open sets. Chen and Tokle~\cite{CT}
derived two-sided global heat kernel estimates for symmetric stable
processes in two classes of unbounded $C^{1,1}$ open sets. See
\cite{BGR} for Varopoulos-type two-sided heat kernel estimates for
symmetric stable processes in a general class of domains including
Lipschitz domains expressed in terms of the surviving probability
function $\P_x(\tau_D>t)$.

This paper can be viewed as a natural continuation of our previous
works~\cite{CKS,CKS1}. We point out that, although this paper adopts
its main strategy from~\cite{CKS}, there are many new difficulties
and differences between obtaining estimates on the transition
densities of relativistic stable processes in open sets and those of
symmetric stable processes and censored stable processes in open
sets studied~\cite{CKS,CKS1}. For example, unlike symmetric stable
processes and censored stable processes, relativistic stable
processes do not have the scaling property, which is one of the main
ingredients used in the approaches of~\cite{CKS,CKS1}. As in
\cite{CKS,CKS1}, the L\'{e}vy system of $X^m$, which describes how
the process jumps [see \eqref{e:levy}], is the basic tool used
throughout our argument because $X^m$ moves by ``pure jumps.''
However, the L\'{e}vy density of $X^m$ does not have a simple form and
has exponential decay at infinity as opposed to the polynomial decay
of the L\'{e}vy density of symmetric stable processes. [See
\eqref{e:jm}--\eqref{e:2.10} and \eqref{H:1}--\eqref{H:2} below.]
Moreover, in this paper we aim at obtaining sharp estimates that
are uniform in $m \in(0, M]$; that is, the constants $C_1$ and
$C_2$ in Theorem~\ref{t:main} are independent of $m \in(0, M]$.
This requires very careful and detailed estimates throughout our
proofs.

The approach of this paper uses a combination of probabilistic and
analytic techniques, but it is mainly probabilistic. It was first
established in~\cite{R}, and then in~\cite{CS4} by using a different
method, that the Green function of $X^m$ in a bounded $C^{1,1}$ open
set $D$ is comparable to that of $X$ in $D$. We show in Theorem
\ref{t:ufgfnest} below, following the approach of~\cite{CS4}, that
such a~comparison is uniform in $m\in(0, M]$ for small balls. This
uniform Green function estimate is then used to get the boundary
decay rate of $p^m_D(t, x, y)$. When~$x$ and $y$ are far from the
boundary in a scale given by $t$, the near diagonal lower bound
estimate of $p^m_D(t, x,y)$ is derived from the uniform parabolic
inequality (Theorem~\ref{T:2.3}), the uniform exit time estimate
(Theorem~\ref{T:2.2}) and the fact that~$X^m_t$ moves from $x$ to a
neighborhood of $y$ by one single jump with positive probability.
These estimates can be used to get the lower bound estimate on the
global heat kernel $p^m(t, x, y)$. The upper bound estimate on
$p^m(t, x, y)$ is obtained from the heat kernel of Brownian motion
through subordination. This sharp two-sided estimates on the
transition density $p^m(t, x, y)$ in bounded time interval are
presented in Theorem~\ref{T:lbRd} and will be used to derive upper
bound estimates on $p^m_D(t, x, y)$. The estimates in Theorem
\ref{T:lbRd} sharpen the corresponding estimates established earlier
in~\cite{CK2} that are applicable for more general jump processes
with exponentially decaying jump kernels. After the first version of
this paper was written and posted on the arXiv, the authors were
informed that the estimates in Theorem~\ref{T:lbRd} are also
obtained in~\cite{Sz2}. Since $X^m$ can be obtained from $X$ by
pruning jumps in a suitable way (see~\cite{BBCK}, Remarks 3.4 and~3.5), we can conclude that $p^m_D(t, x, y) \leq e^{Mt}
p_D(t, x, y)$ for all $m\in(0, M]$. The upper bound estimate on
$p^m_D(t, x, y)$ (Theorem~\ref{t:ub}) is then obtained by using the
L\'{e}vy system formula, comparison with the heat kernel estimates on
exterior balls (Lemma~\ref{L:2.1}), the estimates on $p_D(t, x, y)$
from~\cite{CKS} and the two-sided estimates on $p^m(t, x, y)$
(Theorem~\ref{T:lbRd}).

When $D$ is a bounded $C^{1,1}$ open set, integrating the
estimates on $p^m_D(t, x, y)$ from Theorem~\ref{t:main} over $t$
yields sharp two-sided sharp estimates on the Green function
$ G^m_D(x, y)
:=\int_0^\infty p^m_D(t, x, y)\,dt$. To state this result, we define a
function $V^{\alpha}_D$ on $D\times D$ by
%
%
\begin{equation}\label{e:1.4}
V^{\alpha}_D(x, y) :=
\cases{
\displaystyle \biggl(1\wedge\frac{
\delta_D(x)^{\alpha/2} \delta_D(y)^{\alpha/2}}{ |x-y|^{\alpha}}
\biggr) |x-y|^{\alpha-d} ,\vspace*{2pt}\cr\qquad\mbox{when }
d>\alpha,
\vspace*{5pt}\cr\displaystyle
\log\biggl( 1+ \frac{ \delta_D(x)^{1/2} \delta_D
(y)^{1/2}}{|x-y|} \biggr) ,\vspace*{2pt}\cr\qquad\mbox{when }
d=1=\alpha, \vspace*{2pt}\cr\displaystyle
( \delta_D(x) \delta_D (y) )^{(\alpha-1)/2}
\wedge\frac{ \delta_D(x)^{\alpha/2} \delta_D (y)^{\alpha/2}}{
|x-y|} ,\cr\qquad\mbox{when } d=1<\alpha.
}
\end{equation}

\begin{theo}\label{T:1.2}
Let $M>0$ be a constant and $D$ a bounded $C^{1,1}$ open set in
$\bR^d$. Then there exists a constant $C_3>1$ depending only on $d,
\alpha, r_0, \Lambda_0$, $M, T$, $\operatorname{diam}(D)$ such that
for every
$m\in(0, M]$ and $(x, y) \in D\times D$,
\[
C_3^{-1} V^{\alpha}_D(x, y)\le G^m_D(x, y)\le C_3 V^{\alpha}_D(x, y).
\]
\end{theo}

The proof of Theorem~\ref{T:1.2} is the same as that of
\cite{CKS}, Corollary 1.2. Theorem~\ref{T:1.2} extends and improves
the Green function estimates obtained in~\cite{CS4,KL,R} in the
sense that our estimates are uniform in $m\in(0, M]$ and the case
$d=1$ is now covered.
Although we do not yet have large time heat kernel estimates when $D$ is
unbounded, the short time heat kernel estimates in Theorem
\ref{t:main}(i) can be used together with
the two-sided Green function estimates on the upper
half space from~\cite{GR2} and a comparison idea from~\cite{CT}
to obtain sharp two-sided estimates on the Green function
$G^m_D(x, y)$ when $D$ is a half-space-like $C^{1, 1}$ open set.
We will address this in a separate paper~\cite{CKS3}.

The rest of the paper is organized as follows. In Section
\ref{S:int} we recall some basic facts about the relativistic stable
process $X^m$ and prove some preliminary uniform results on $X^m$
including
uniform estimates
on the Green function $G^m_D$ for small balls
and annuli,
and uniform parabolic
Harnack inequality.
Some preliminary lower bound of $p^m_D(t, x, y)$ is proved in
Section~\ref{S:ub}, while the proof of Theorem~\ref{t:main} is given in Section~\ref{S:lb}.

In the remainder of this paper, we assume that $m>0$. We will use
capital letters $C_1, C_2, \ldots$ to denote constants in the
statements of results, and their labeling will be fixed. The lower
case constants $c_1, c_2, \ldots$ will denote generic constants used
in proofs, whose exact values are not important and can change from
one appearance to another. The labeling of the lower case constants
starts anew in every proof.
The dependence of the lower case constants on the dimension $d$ will
not be mentioned explicitly. We will use $\partial$ to denote a
cemetery point and for every function $f$, we extend its definition
to~$\partial$ by setting $f(\partial)=0$. We will use $dx$ to
denote the Lebesgue measure in~$\bR^d$. For a Borel set $A\subset
\bR^d$, we also use $|A|$ to denote its Lebesgue measure.\looseness=-1

\section{Relativistic stable processes and
some uniform estimates}\label{S:int}

The L\'{e}vy measure of the relativistic $\alpha$-stable process
$X^m$, defined in \eqref{e:ch}, has a density
%
%
\begin{eqnarray}\label{e:jm}
J^m(x)&\hspace*{3pt}=&
j^m(|x|)\nonumber
\\[-8pt]
\\[-8pt]&:=&\frac{\alpha}{2 \Gamma(1-\fraca{\alpha}2)} \int
_0^{\infty}
(4\pi u)^{-d/2}e^{-\fraca{|x|^2}{4u}}e^{-m^{2/\alpha}u}
u^{-(1+\fraca{\alpha}{2})}\,du,
\nonumber
\end{eqnarray}
which is continuous and radially decreasing on $\bR^d\setminus\{0\}$
(see~\cite{R}, Lemma~2).
Here and in the rest of this paper, $\Gamma$ is the Gamma function
defined by $\Gamma(\lambda):= \int^{\infty}_0 t^{\lambda-1}
e^{-t}\,dt$ for every $\lambda> 0$.
Put $J^m(x,y):= j^m(|x-y|)$. Let\break $ {\cal A}(d, -\alpha):=
\alpha2^{\alpha-1}\pi^{-d/2}
\Gamma(\frac{d+\alpha}2) \Gamma(1-\frac{\alpha}2)^{-1}$.\vadjust{\goodbreak}
Using change of variables twice, first with $u=|x|^2v$ then with
$v=1/s$, we get
%
%
\begin{equation}\label{e:jm2}
J^m(x,y)= {\cal A} (d, -\alpha) |x-y|^{-d-\alpha} \psi
(m^{1/\alpha}|x-y|),
\end{equation}
where
%
%
\begin{equation}\label{e:psi} \psi(r):= 2^{-(d+\alpha)} \Gamma
\biggl(
\frac{d+\alpha}{2} \biggr)^{-1} \int_0^\infty s^{\fracb{d+\alpha
}{2}-1} e^{-\fraca{s}{ 4} -\fraca{r^2}{ s} } \, ds,
\end{equation}
which satisfies
$\psi(0)= 1$ and
%
%
\begin{equation}\label{e:2.10}
c_1^{-1}e^{-r}r^{(d+\alpha-1)/2} \le\psi(r) \le c_1
e^{-r}r^{(d+\alpha-1)/2} \qquad\mbox{on } [1, \infty)
\end{equation}
for some $c_1>1$ (see~\cite{CS4}, pages 276--277, for details).
In particular, we see that for $m>0$,
%
%
\begin{equation}\label{e:2.5}
J^m(x, y)= m^{(d+\alpha)/\alpha} J^1(m^{1/\alpha}x, m^{1/\alpha} y).
\end{equation}
We denote the L\'{e}vy density of $X$ by
\[
J(x, y):= J^0(x,y)={\cal A}(d, -\alpha)|x-y|^{-(d+\alpha)}.
\]

The L\'{e}vy density gives rise to a L\'{e}vy system for $X^m$, which
describes the jumps of the process $X^m$:
for any $x\in\bR^d$, stopping time $T$ (with respect to the
filtration of $X^m$) and nonnegative measurable function $f$ on
$\bR_+ \times\bR^d\times\bR^d$ with $f(s, y, y)=0$ for all $y\in
\bR^d$ and $s\ge0$,
%
%
\begin{equation}\label{e:levy}
 \E_x \biggl[\sum_{s\le T} f(s,X^m_{s-}, X^m_s)
\biggr] = \E_x \biggl[ \int_0^T \biggl(
\int_{\bR^d} f(s,X^m_s, y) J^m(X^m_s,y)\,dy \biggr)\,ds \biggr].\hspace*{-35pt}
\end{equation}
(See, e.g.,~\cite{CK}, proof of Lemma 4.7,
and~\cite{CK2}, Appendix~A.)

For $r\in(0, \infty)$, we define
%
%
\begin{equation}\label{e:xi}
\xi(r):=
\cases{\displaystyle
r^2,&\quad when $d+\alpha>2$, \vspace*{2pt}\cr\displaystyle
r^{1+\alpha}
,&\quad when $d=1>\alpha$, \vspace*{2pt}\cr\displaystyle
r^2
\ln\biggl(1+\frac{1}{r}\biggr)
,&\quad when $d=1=\alpha$.
}
\end{equation}

We start with an elementary inequality.

\begin{lemma}\label{l:bdonpsi}
For any $R_0>0$, there exists $C_4=C_4(d, \alpha, R_0)>0$
such that for all $r\in(0, R_0]$,
\[
1-\psi(r)\le C_4 \xi(r).
\]
\end{lemma}

\begin{pf}
We have
\[
1-\psi(r)=2^{-(d+\alpha)} \Gamma\biggl( \frac{d+\alpha}{2}
\biggr)^{-1} \biggl(\int_0^{r^2}
+\int_{r^2}^\infty\biggr)
s^{\fracb{d+\alpha}{ 2}-1} e^{-\fraca{s}{ 4}}(1- e^{-\fraca{r^2}{
s} })
\, ds.
\]
Note that
%
%
\begin{equation}\label{e:psi1}\qquad
\int_0^{r^2}s^{\fracb{d+\alpha}{ 2}-1} e^{-\fraca{s}{ 4}}(1 -
e^{-\fraca{r^2}{ s} })\, ds\le\int_0^{r^2} s^{\fracb{d+\alpha
}{2}-1}\,ds
\le c_1r^{d+\alpha},
\end{equation}
and that, by
the inequality $1-e^{-z} \le z$ for $z \ge0$,
%
%
\begin{eqnarray}\label{e:psi2}
\int_{r^2}^\infty s^{\fracb{d+\alpha}{ 2}-1} e^{-\fraca{s}{ 4}}(1-
e^{-\fraca{r^2}{ s} }) \, ds&\le& r^2 \int_{r^2}^\infty s^{\fracb
{d+\alpha}{
2}-2} e^{-\fraca{s}{ 4}} \, ds\nonumber
\\[-8pt]
\\[-8pt]&\le& c_2 \xi(r).
\nonumber
\end{eqnarray}
We arrive at the conclusion of this lemma by combining
\eqref{e:psi1} and \eqref{e:psi2}.
\end{pf}

The next two inequalities, which can be seen easily from
the monotonicity of $\psi$ and \eqref{e:2.10}, will be used several
times in this paper. For any $a>0$ and $M>0$, there exist positive
constants $C_{5}$ and $C_{6}$ depending only on~$a$ and~$M$ such
that for any $m\in(0, M]$,
%
%
\begin{equation}\label{H:1}
j^m(r)\le C_{5}j^m(2r) \qquad\mbox{for every } r\in
(0,a]
\end{equation}
and
%
%
\begin{equation}\label{H:2}
j^m(r)\le C_{6} j^m(r+a) \qquad\mbox{for every } r>a.
\end{equation}

We will use $p^m(t, x, y)=p^m(t, x-y)$ to denote the transition
density of~$X^m$ and use $p(t, x, y)$ to denote the
transition density of $X$. It is well known
that (cf.~\cite{CK})
%
%
\begin{equation}\label{e:1.1}
p (t, x, y) \asymp t^{-d/\alpha} \wedge
\frac{t}{|x-y|^{d+\alpha}} \qquad\mbox{on } (0,
\infty)\times\bR^d \times\bR^d.
\end{equation}

Here and in the sequel,
for two nonnegative functions $f, g$, $f\asymp g$
means that there is a positive constant $c_0>1$ so that
$c_0^{-1} f\leq g\leq c_0f$ on their common domain of definitions.
It is also known that
%
%
\begin{equation}\label{e:m11}
p^1 (t, x) = e^t\int_0^\infty(4\pi
u)^{-d/2}e^{- |x|^2/(4u)}e^{-u} \theta_\alpha(t,u)\,du,
\end{equation}
where $\theta_\alpha(t,u)$ is the transition density of an
$\frac{\alpha}{2}$-stable subordinator with the Laplace transform
$e^{-t \lambda^{\alpha/2}}$. It follows from
\cite{BG}, Theorem 2.1, and~\cite{Z}, (2.5.17), (2.5.18), that
\[
\theta_\alpha(t,u) \le c t u^{-1-\alpha/2} \qquad\mbox{for every }
t
>0, u >0.
\]
Thus by \eqref{e:jm} and \eqref{e:m11}, there exists $L=L(\alpha)>0$
such that
%
%
\begin{equation}\label{e:ma1}
\qquad p^1 (t, x, y) \le L {t e^t }J^{1}(x,y) \qquad\mbox{for all }
(t,x,y) \in(0, \infty)\times\bR^d \times\bR^d.
\end{equation}
From \eqref{e:ch}, one can easily see that $X^m$ has the following
approximate scaling property:
$\{m^{-1/\alpha}(X^1_{m t}-X^1_0), t\geq
0\}$ has the same distribution as that of $ \{X^m_t-X^m_0,
t\geq0 \}.$
In terms of transition densities, this approximate scaling property
can be written as
%
%
\begin{equation}\label{scale_p}
p^m(t,x,y) = m^{d/\alpha}p^1 (mt, m^{1/\alpha} x, m^{1/\alpha} y).
\end{equation}
Thus by \eqref{e:2.5},
\eqref{e:ma1} and \eqref{scale_p},
we have
%
%
\begin{equation}\label{e:ma2}
\qquad p^m (t, x, y) \le L t e^{mt} J^{m}(x,y)
\qquad\mbox{for }
(t,x,y) \in(0, \infty)\times\bR^d \times\bR^d .
\end{equation}
On the other hand, by~\cite{R}, Lemma 3,
there exists
$c=c(\alpha)>0$ such that
%
%
\begin{equation}\label{e:mb2}
p^m (t, x, y) \le c (m^{d/\alpha-d/2}t^{-d/2}+t^{-d/\alpha}).
\end{equation}

For any open set $D$, we use $\tau^m_D$ to denote the first
exit time from $D$ for $X^m$, that is, $\tau^m_D=\inf\{t>0\dvtx
X^m_t\notin
D\}$ and let $\tau_D$ be the first exit time from $D$ for~$X$.
We define $X^{m,D}$ by $X^{m, D}_t(\omega)=X^m_t(\omega)$ if $t<
\tau^m_D(\omega)$ and $X^{m, D}_t(\omega)=\partial$ if $t\geq
\tau^m_D(\omega)$. We define $X^D$ similarly. $X^{m,D}$ is called
the subprocess of $X^m$ killed upon exiting $D$ (or, the killed
relativistic stable process in $D$ with mass $m$), and $X^D$ is
called the killed symmetric $\alpha$-stable process in $D$.

It is known (see~\cite{CK2}) that $X^{m,D}$ has a transition
density $p^m_D(t, x, y)$, which is continuous on $(0, \infty) \times
D \times D$
with respect to the Lebesgue measure. Note that
the transition density $p^m_D(t, x, y)$ may not be continuous on
$\overline{ D} \times\overline{ D} $ if the boundary of $D$ is
irregular.

We will use $G^m_D(x, y):=\int_0^{\infty} p^m_D(t, x, y)\,dt$ to
denote the Green function of~$X^{m,D}$.
We use $p_D(t, x, y)$ and $G_D(x, y)$ to denote the transition
density
and the Green function of $X^D$, respectively.

The Dirichlet heat kernel
$p^m_D(t, x, y)$ also has the following approximate scaling property:
%
%
\begin{equation}\label{scale_kp}
p_D^m(t,x,y) = m^{d/\alpha}p_{m^{1/\alpha}D}^1 (mt, m^{1/\alpha}
x, m^{1/\alpha} y).
\end{equation}
Thus the Green function $G^m_D(x, y)$ of $X^{m,D}$ satisfies
%
%
\begin{equation}\label{scale_kg}
 G^m_{D}(x, y) =
m^{(d-\alpha)/\alpha} G^1_{m^{1/\alpha}D} (m^{1/\alpha} x,
m^{1/\alpha} y) \qquad\mbox{for every } x, y\in D.\hspace*{-35pt}
\end{equation}

\begin{remark}\label{R:2.2}
We point out here that the uniform heat kernel estimates in Theorem
\ref{t:main} do not follow from a combination of the sharp heat kernel
estimates of $p^1_D(t, x, y)$ and the scaling property
\eqref{scale_kp}. This is because if
$D$ is a $C^{1,1}$ open sets with $C^{1,1}$ characteristics $(r_0,
\Lambda_0)$, then $m^{1/\alpha}D$ is a $C^{1,1}$ open sets with
different $C^{1,1}$ characteristics $(m^{1/\alpha}r_0, m^{-1/\alpha
}\Lambda_0)$.
\end{remark}

Let
\[
J_m(x,y)=J(x,y)-J^m(x,y) = {\cal A} (d, -\alpha)
|x-y|^{-d-\alpha} \bigl(1-\psi(m^{1/\alpha}|x-y|)\bigr).
\]
Then
%
%
\begin{equation}\label{e:dfjm}
\int_{\bR^d} J_m(x,y)\,dy=m \qquad\mbox{for all } x\in\bR^d.
\end{equation}
(See~\cite{R}, Lemma 2.) Thus $X^m$ can be constructed from $X$ by
reducing jumps via Meyer's construction (see~\cite{BBCK}, Remarks 3.4
and 3.5). By~\cite{R}, Lemma 5, or
\cite{BBCK}, (3.18), we have
%
%
\begin{equation}\label{e:Meyer}
 p^m_D(t,x,y) \le e^{mt} p_D(t,x,y)
\qquad\mbox{for every } (t,x,y) \in(0, \infty)\times D \times D.\hspace*{-35pt}
\end{equation}

In the next two results, we discuss the Green function of
one-dimensional symmetric
$\alpha$-stable processes killed upon exiting $B=(0, 2)\subset\bR$.
Define for $x, y\in B$,
%
%
\begin{equation}\label{e:f1}
f(x, y):= \frac{\delta_B(x) \delta_B(y)}{|x-y|^2}.
\end{equation}

\begin{lemma}\label{l:g3g1}
Suppose that
$B=(0, 2) \subset\bR$
and $\alpha>1$.

\begin{longlist}[(ii)]
\item[(i)]
There exists $C_{7}=C_{7}(\alpha)>0$ such that
\[
\frac{G_B(x,y) G_{B}(y,z)} { G_{B}(x,z)} \le C_{7}
\qquad\mbox{for every } x,y,z \in B.
\]
\item[(ii)]
If $f(x, w)\ge4$, there exists $C_{8}=C_{8}(\alpha)>0$ such that
\[
\frac{G_B(x,y) G_{B}(z,w)} { G_{B}(x,w)} \le C_{8}
\delta_B (y)^{(\alpha-1)/2} \delta_B
(z)^{(\alpha-1)/2} \le C_{8} \qquad\mbox{for } y, z\in
B.
\]
\end{longlist}
\end{lemma}

\begin{pf}
(i) follows from~\cite{BB}, (3.5). So we only need to prove (ii).

Note that (see~\cite{CZ}, page 187)
$|x-y| \le\delta_B(x) \wedge\delta_B(y)$ and
$\delta_B(x) \wedge\delta_B(y)\ge\frac12 (\delta_B(x) \vee
\delta_B(y) )$ if $ f(x,y) \ge4$.
We know from
\cite{BB}, Corollary 3.2, or~\cite{CKS}, Corollary 1.2,
that
%
%
\begin{equation}\label{e:Gb} G_B(x,
y) \asymp( \delta_B(x) \delta_B (y) )^{(\alpha-1)/2}
\wedge
\frac{ \delta_B(x)^{\alpha/2} \delta_B (y)^{\alpha/2}}{ |x-y|}.
\end{equation}
So when $f(x, w)\geq4$, we have by \eqref{e:Gb} that
\begin{eqnarray*}
\frac{G_B(x,y) G_{B}(z,w)} { G_{B}(x,w)} &\le&
c_1 \frac{ ( \delta_B(x) \delta_B (y) )^{(\alpha-1)/2}
( \delta_B(z) \delta_B (w) )^{(\alpha-1)/2}}
{ (\delta_B(x) \delta_B (w))^{(\alpha-1)/2}} \\
&=& c_1
\delta_B (y)^{(\alpha-1)/2}\delta_B (z)^{(\alpha-1)/2} .
\end{eqnarray*}
\upqed
\end{pf}

The second part of the next result strengthens~\cite{BB}, (3.4).

\begin{lemma}\label{l:g3g2}
Suppose that
$B=(0, 2) \subset\bR$ and $\alpha=1$. Let $f$ be as in~\eqref{e:f1} and define
$
F(x,y):= \log( 1+ f(x,y)^{1/2} ).
$
\begin{longlist}[(ii)]
\item[(i)]
If $f(x,w) \ge4$, there exists $C_{9}>0$ such that
\[
\frac{G_B(x,y) G_{B}(z,w)} { G_{B}(x,w)} \le C_{9} F(x, y)F(z, w),
\qquad y, z\in B.
\]
\item[(ii)] There exists $C_{10}>0$ such that
\[
\frac{G_B(x,y) G_{B}(y,z)} { G_{B}(x,z)} \le C_{10}\bigl(1+ F(x,y) +
F(y,z)\bigr),
\qquad x, y, z\in B.
\]
\end{longlist}
\end{lemma}

\begin{pf}
(i) is an immediate consequence of~\cite{CKS}, Corollary 1.2. Using~\cite{CKS},
Corollary 1.2, (ii) can be proved by following the
argument of the proof of~\cite{CZ}, Theorem~6.24. We omit the
details.
\end{pf}

For $r\in(0, 1]$, we define
\[
\sigma(r)=
\cases{\displaystyle
r^{2-\alpha-d} ,&\quad when $d+\alpha>2$, \cr\displaystyle
1
,&\quad when $d=1>\alpha$, \cr\displaystyle
\ln(1+1/r)
,&\quad when $d=1=\alpha$.
}
\]
The following result will be used to prove Theorem~\ref{t:ufgfnest}.
Note that the case $d=1\le\alpha$ in Lemma~\ref{l:integral}(i) does
not follow from~\cite{GR}, Lemma 3.14.
%
\begin{lemma}\label{l:integral}
\textup{(i)} If $B$ is a ball of radius $1$ in $\bR^d$, then,
\[
\sup_{x, y\in B, x\neq y}\int_{B\times B} \frac{G_{B}(x,
w)\sigma(|w-z|)G_{B}(z, y)}{G_{B}(x, y)}\,dw\,dz<\infty.
\]

\textup{(ii)} If $d\ge2$ and $U$ is an annulus of inner radius $1$ and outer
radius $3/2$ in $\bR^d$, then
\[
\sup_{x, y\in U, x\neq y}\int_{U\times U} \frac{G_{U}(x,
w)\sigma(|w-z|)G_{U}(z, y)}{G_{U}(x, y)}\,dw\,dz<\infty.
\]
\end{lemma}

\begin{pf} We only present the proof of (i). The proof of (ii) is similar
to the proof of (i) for the case $d>\alpha$. We prove (i) by dealing
with two separate cases.

Case 1: $d>\alpha$. In this case, by repeating the argument in
\cite{CS7}, Example 2 (also see~\cite{GR}, Lemma 3.14),
we know that there exists $c_1=c_1(d,
\alpha)>0$ such that
\begin{eqnarray*}
&&\frac{G_{B}(x, w)\sigma(|w-z|)G_{B}(z, y)}{G_{B}(x, y)}\\
&& \qquad \le c_4 \biggl(\frac1{|z-y|^{d-\alpha/2}|w-z|^{d+\alpha
-\beta}}+
\frac1{|x-w|^{d-\alpha/2}|w-z|^{d+\alpha-\beta}} \\
&& \hphantom{c_4 \biggl(}\qquad \quad{}+
\frac1{|z-y|^{d-\alpha}|w-z|^{d+\alpha-\beta}} +
\frac1{|x-w|^{d-\alpha}|w-z|^{d+\alpha-\beta}} \\
&& \hphantom{c_4 \biggl(}\qquad \quad{} + \frac1{|x-w|^{d-\alpha/2}|z-y|^{d-\alpha
/2}|w-z|^{3\alpha
/2-\beta}}
\\
&&\hspace*{96pt} \qquad \quad{}
+
\frac1{|x-w|^{d-\alpha/2}|z-y|^{d-\alpha}|w-z|^{2\alpha-\beta
}} \biggr),
\end{eqnarray*}
where $\beta=2$ when $d\ge2$ and $\beta=1+\alpha$ when
$d=1>\alpha$. The conclusion now follows immediately.

Case 2: $d=1\le\alpha$. In this case, it follows from the first
part of the proof of
\cite{GR}, Proposition 3.17, that
\[
\sup_{x, y\in B, x\neq y, f(x, y)\le4}\int_{B\times B}
\frac{G_{B}(x, w)\sigma(|w-z|)G_{B}
(z, y)}
{G_{B}(x, y)}\,dw\,dz <
\infty,
\]
where the $f$ is the function defined in \eqref{e:f1}. The
inequality
\[
\sup_{x, y\in B, x\neq y, f(x, y)\ge4}\int_{B\times B}
\frac{G_{B}(x, w)\sigma(|w-z|)G_{B}
(z, y)}{G_{B}(x, y)}\,dw\,dz<\infty
\]
follows easily from Lemmas~\ref{l:g3g1} and~\ref{l:g3g2}.
\end{pf}

The following result will be used later in this paper.
Note that this result does not follow from the main result in
\cite{GR}, since the constants in the following results are uniform
in $m\in(0, \infty)$ and $r\in(0, R_0 m^{-1/\alpha}]$.
It is known that (see~\cite{CS1} and~\cite{K}
for the case $d\ge2$
and~\cite{BB} for the case $d=1$) that $G_B(x, y)$ is comparable to
$V^\alpha_B(x, y)$ of \eqref{e:1.4}.

\begin{theo}\label{t:ufgfnest}
There exist positive constants
$R_0\in(0, 1]$
and
$C_{11}>1$ depending only on $d$ and $\alpha$ such that for any
$m\in(0, \infty)$, any ball $B$ of radius $r\le R_0 m^{-1/\alpha}$,
\[
C_{11}^{-1}G_B(x, y)\le G^m_B(x, y)\le C_{11}G_B(x, y), \qquad x, y\in B.
\]
Furthermore, in the case $d\ge2$, there exists a constant
$C_{12}=C_{12}(d, \alpha)>1$ such that for any $m\in(0, \infty)$,
$r\in(0, R_0 m^{-1/\alpha}]$ and any annulus $U$ of inner radius
$r$ and outer radius $3r/2$,
\[
C_{12}^{-1}G_U(x, y)\le G^m_U(x, y)\le
C_{12}G_U(x, y), \qquad x, y\in U.
\]
\end{theo}

\begin{pf}
We only present the proof for balls, the case of annuli is similar.
By~\cite{BB,CKS,K}, $G_B(x, y) \asymp V^\alpha_B (x, y)$. Hence by
\eqref{scale_kg}, we only need to prove the theorem for $m=1$. In
this proof we will use $B_r$ to denote the ball $B(0, r)$.

Put
\[
F(x, y):=\frac{J^1(x, y)}{J(x, y)}-1=\psi(|x-y|)-1
, \qquad x, y\in\bR^d.
\]
Then it follows from \eqref{e:jm}--\eqref{e:psi} that there exists
$c_1=c_1(d, \alpha)>0$ such that for
any $r\in(0, 1]$, $ \inf_{x, y\in B_r}F(x, y)\ge c_1-1$. It
follows from Lemma~\ref{l:bdonpsi} that there exists $c_2=c_2( d,
\alpha)>0$ such that for any $r\in(0, 1]$ and $x, y\in B_{1}$,
%
%
\begin{eqnarray}\label{e:bdonf}
  &&|F(rx, ry)|+\bigl|\ln\bigl(1+F(rx, ry)\bigr)\bigr|+ \bigl(e^{4|\ln(1+F(rx, ry))|}-1\bigr)
  \nonumber
  \\[-8pt]
  \\[-8pt]
  && \qquad \le c_2
\xi(r|x-y|).
\nonumber
\end{eqnarray}
For $ x\in B_r$, put
\[
q_{B_r}(x):=
\int_{B_r^c}J_1(x,y)\,dy = {\cal A} (d, -\alpha)
\int_{B_r^c}|x-y|^{-d-\alpha} \bigl(1-\psi
(|x-y|)\bigr)\,dy.
\]

Then it follows from~\cite{CS4}, Section 3, that
\[
G^1_{B_r}(x, y)=G_{B_r}(x, y)\E_x^y[K^{B_r}(\tau_{B_r})]
\qquad\mbox{for every } x, y\in B_r,
\]
where
\begin{eqnarray*}
K^{B_r}(t)&:=&\exp\biggl(\sum_{0<s\le t}\ln\bigl(1+F(X^{B_r}_{s-},
X^{B_r}_s)\bigr)
\\&&\hphantom{\exp\biggl(\sum_{0<s\le t}}{}-\int^t_0\int_{B_r}F(X^{B_r}_s, y)J(X^{B_r}_s,
y)\,dy\,ds + \int^t_0q_{B_r}( X^{B_r}_s)\,ds \biggr).
\end{eqnarray*}

Using the scaling property of $G_{B_r}$, we get
%
%
\begin{eqnarray} \label{aa}
&&  \sup_{x, y\in B_r, x\neq y}\int_{B_r\times B_r}
\frac{G_{B_r}(x,
w)(e^{4|\ln(1+F(w, z))|}-1)G_{B_r}(z, y)}{G_{B_r}(x,
y)|w-z|^{d+\alpha}}\,dw\,dz\nonumber\hspace*{-35pt}
\\[-4pt]
\\[-12pt]
&&  \qquad=
\sup_{x, y\in B_{1}, x\neq y}\int_{B_{1}\times B_{1}} \frac{G_{B_{1}}(x,
w)(e^{4|\ln(1+F(rw, rz))|}-1)G_{B_1}(z, y)}{G_{B_1}(x,
y)|w-z|^{d+\alpha}}\,dw\,dz,
\nonumber\hspace*{-35pt}
\end{eqnarray}
%
%
\begin{eqnarray}
\label{bb}
&&\sup_{x, y\in B_r, x\neq y}\int_{B_r\times B_r} \frac{G_{B_r}(x,
w)|F(w, z)|G_{B_r}(z, y)}{G_{B_r}(x,
y)|w-z|^{d+\alpha}}\,dw\,dz \nonumber
\\[-8pt]
\\[-8pt]
&& \qquad=\sup_{x, y\in B_1, x\neq y}\int_{B_1\times B_1} \frac{G_{B_1}(x,
w)|F(rw, rz)|G_{B_1}(z, y)}{G_{B_1}(x,
y)|w-z|^{d+\alpha}}\,dw\,dz
\nonumber
\end{eqnarray}
and
%
%
\begin{eqnarray}\label{cc}
&&\sup_{x, y\in B_r, x\neq y}\int_{B_r}\frac{G_{B_r}(x, w)G_{B_r}(w,
y)}{G_{B_r}(x, y)}q_{B_r}(w)\,dw\nonumber
\\[-8pt]
\\[-8pt]
&& \qquad=r^{\alpha}\cdot\sup_{x, y\in B_1,
x\neq y}\int_{B_1}\frac{G_{B_1}(x, w)G_{B_1}(w, y)}{G_{B_1}(x,
y)}q_{B_r}(rw)\,dw.
\nonumber
\end{eqnarray}
Using \eqref{e:bdonf}--\eqref{bb} and Lemma
\ref{l:integral}, we have for $r \in(0,1]$,
\begin{eqnarray*}
\sup_{x, y\in B_r, x\neq y}\int_{B_r\times B_r} \frac{G_{B_r}(x,
w)(e^{4|\ln(1+F(w, z))|}-1)G_{B_r}(z, y)}{G_{B_r}(x,
y)|w-z|^{d+\alpha}}\,dw\,dz \le c_3 r
\end{eqnarray*}
and
\[
\sup_{x, y\in B_r, x\neq y}\int_{B_r\times B_r} \frac{G_{B_r}(x,
w)|F(w, z)|G_{B_r}(z, y)}{G_{B_r}(x,
y)|w-z|^{d+\alpha}}\,dw\,dz\le c_3 r .
\]
By applying \eqref{e:dfjm}, the 3G inequality [Lemma~\ref{l:g3g1}(ii)
and Lemma~\ref{l:g3g2}(ii) for $d=1$] and
\eqref{cc}, we
also have
\[
\sup_{x, y\in B_r, x\neq y}\int_{B_r}\frac{G_{B_r}(x,
w)G_{B_r}(w, y)}{G_{B_r}(x, y)}q_{B_r}(w)\,dw \le c_3r^\alpha.
\]
Now choose $R_0>0$ small enough so that for $r \le R_0$,
\begin{eqnarray*}
&\displaystyle\sup_{x, y\in B_r, x\neq y}\int_{B_r\times B_r} \frac
{G_{B_r}(x,
w)(e^{4|\ln(1+F(w, z))|}-1)G_{B_r}(z, y)}{G_{B_r}(x,
y)|w-z|^{d+\alpha}}\,dw\,dz\le\frac12,&
\\
&\displaystyle\sup_{x, y\in B_r, x\neq y}\int_{B_r\times B_r} \frac
{G_{B_r}(x,
w)|F(w, z)|G_{B_r}(z, y)}{G_{B_r}(x,
y)|w-z|^{d+\alpha}}\,dw\,dz\le\frac18&
\end{eqnarray*}
and
\[
\sup_{x, y\in B_r, x\neq y}\int_{B_r}\frac{G_{B_r}(x,
w)G_{B_r}(w, y)}{G_{B_r}(x, y)}\,dw\le\frac18.
\]
Using the three displays above, we can follow the argument in
\cite{CS7}, Proposition~2.3
(with the constants involved there taken to be $\alpha=\gamma=2,
\theta=1/2$)
to conclude that for all $r \le R_0$,
\[
\sup_{x, y\in B_r, x\neq y}\E^y_x [K^{B_r}(\tau_{B_r})
]\le
2^{3/4}.
\]
Now the upper bound on $G^1_{B_r}$ follows immediately. The lower
bound on~$G^1_{B_r}$ is an easy consequence of Jensen's inequality
(see~\cite{CS7}, Remark 2, for details).
\end{pf}

In the remainder of this paper,
$R_0\in(0, 1]$ will always stand for the
constant in Theorem~\ref{t:ufgfnest}.
The next corollary will be used in Section~\ref{S:lb}.

\begin{cor}\label{c:ufgfnest1}
There exist positive constants $C_{13}>1$ and $C_{14}<1$ depending
only on $d$ and $\alpha$ such that for any $m\in(0, \infty)$,
any $r\le R_0 m^{-1/\alpha}$, any ball $B$ of radius $r$ and,
when $d\ge2$, any annulus $U=B(x_0, 3r/2) \setminus\overline{B(x_0, r)}$
%
%
\begin{eqnarray}\label{ce1}
\P_x ( X^m_{\tau^m_B} \in A ) &\le& C_{13}
\P_x ( X_{\tau_B} \in A )
\qquad\mbox{for every } x \in B \mbox{ and } A \subset B^c,\hspace*{-35pt}
\\
\label{ce2}
\P_x ( X^m_{\tau^m_U} \in A ) &\le& C_{13}
\P_x ( X_{\tau_U} \in A )
\qquad\mbox{for every } x \in U \mbox{ and } A \subset U^c.\hspace*{-35pt}
\end{eqnarray}
In addition,
if $N \ge2R_0$,
then for every $ x \in U$ and
$A \subset B(x_0, Nm^{-1/\alpha}) \setminus B(x_0, 3r/2)$,
%
%
\begin{equation}\label{ce3}
\P_x ( X^m_{\tau^m_U} \in A ) \ge C_{14} \psi(2R_0+ N) \P_x (
X_{\tau_U} \in A )
.
\end{equation}
\end{cor}

\begin{pf}
By \eqref{e:levy} and~\cite{Sz},
\[
\P_x ( X^m_{\tau^m_B} \in A ) = \int_A \int_B G^m_B(x,y)
J^m(y,z)\,dy\,dz\vadjust{\goodbreak}
\]
and
\[
\P_x ( X^m_{\tau^m_U} \in A ) = \int_A
\int_U G^m_U(x,y) J^m(y,z)\,dy\,dz.
\]
Thus, using Theorem
\ref{t:ufgfnest} and the fact $J^m \le J^0$,
\eqref{ce1} and \eqref{ce2} follow immediately.

Moreover, when $y \in B(x_0, 3r/2) \setminus\overline{B(x_0, r)}$ and
$z \in A \subset
B(x_0, Nm^{-1/\alpha}) \setminus B(x_0, 3r/2)$,
$m^{1/\alpha}|y-z| \le2R_0+ N
.$
Thus
$J^m (y,z) \ge\psi(2R_0+ N) J(y,z)$
and, using Theorem~\ref{t:ufgfnest},
\eqref{ce3} follows.
\end{pf}

Later in this paper, we will also need the following exit time
estimate and parabolic Harnack inequality that are uniform in $m\in
(0, M]$. These results are extensions of Proposition 4.9 and
Theorem~4.12 of~\cite{CK2}, respectively.

\begin{theo}\label{T:2.2}
For any $M>0$, $R>0$, $A> 0$ and $B\in(0, 1)$, there exists $\gamma
=\gamma(A, B, M, R)\in(0, 1/2)$ such that for every $m\in(0,
M]$, $r\in(0, R]$ and $x\in\bR^d$,
\[
\P_x \bigl( \tau^m_{B(x, Ar)}<\gamma r^\alpha\bigr) \le B.
\]
\end{theo}

\begin{pf} Let $Y^m$ be a symmetric pure jump process on $\bR^d$ with jump
kernel given by
\[
J^m_0(x, y) =
\cases{\displaystyle
j^m(|x-y|) ,&\quad if $|x-y|\leq1$, \cr\displaystyle
j^m(1) |x-y|^{-(d+\alpha)} ,&\quad if $|x-y|>1$.
}
\]
Note that $J^m_0(x, y)\geq J^m(x, y)$. In view of
\eqref{e:jm}--\eqref{e:2.10} and \eqref{e:dfjm}, there are constants
$c_i=c_i(M, \alpha)>0, i=1, 2$, such that
%
%
\begin{equation}\label{e:2.26}
\frac{c_2}{|x-y|^{d+\alpha}} \leq J^m_0(x, y)\leq
\frac{c_1}{|x-y|^{d+\alpha}}
\end{equation}
for every $ m\in(0, M]$ and $ x, y\in\bR^d$, and
%
%
\begin{equation}\label{e:2.27}
\sup_{z\in\bR^d} {\cal J}^m_0(z) \leq M \qquad\mbox{for every }
m\in(0, M],
\end{equation}
where ${\cal J}^m_0(z):=\int_{\bR^d} (J^m_0(z,
w)-J^m(z,w) )\,dw$.
In view of \eqref{e:2.26}, it follows from
\cite{CK}, Proposition 4.1, that for each $M>0$, $R>0$, $A>0$ and
$B\in(0, 1)$, there is $\gamma=\gamma(A, B, M, R)\in(0, 1)$
such that for every $m\in(0, M]$, $r\in(0, R]$ and $x\in
\bR^d$,
\[
\P_x \bigl( \tau^{Y^m}_{B(x, Ar)}<\gamma r^\alpha\bigr) \le
B/2,
\]
where $\tau^{Y^m}_{B(x, Ar)}$ is the first time the process $Y^m$
exits the set $B(x, Ar)$. On the other hand, in view of
\eqref{e:2.27}, $Y^m$ can
be obtained from $X^m$
by adding new jumps according to the jump
kernel $J^m_0(x, y)-J^m(x, y)$ through Meyer's construction (see
\cite{BBCK}, Remark 3.4). Hence we have for every $m\in(0, M]$,
$r\in(0, R]$ and $x\in\bR^d$,
\begin{eqnarray*}
&&\P_x \bigl( \tau^m_{B(x, Ar)}<\gamma r^\alpha\bigr) \\
&& \quad\leq\P_x \bigl( \tau^{Y^m}_{B(x, Ar)}<\gamma r^\alpha
\mbox{ and there is no new jumps added to } X^m \mbox{ by time }
\gamma r^\alpha\bigr) \\
&& \quad \quad{} + \P_x
( \mbox{there is at least one new jump added
to }
X^m \mbox{ by time } \gamma r^\alpha) \\
&& \quad\leq B/2 + \bigl( 1-e^{-\gamma r^\alpha \|{\cal
J}^m_0\|_\infty} \bigr) \leq B/2+ ( 1-e^{-\gamma R^\alpha
M} )<B,
\end{eqnarray*}
where the last inequality is achieved by decreasing the value of
$\gamma$ if necessary.
\end{pf}

We now introduce the space--time process $Z^m_s:=(V_s, X^m_s)$, where
$V_s=V_0- s$. The filtration generated by $Z^m$ and satisfying the usual
condition will be denoted as $\{ \widetilde\sF_s; s\geq0\}$.
The law of the space--time process $s\mapsto Z^m_s$ starting from $(t,
x)$ will be denoted as $\mathbb{P}^{(t, x)}$ and as usual,
$\mathbb{E}^{(t, x)}( \cdot) =\int\cdot \mathbb{P}^{(t,
x)} (d \omega).$

We say that a nonnegative Borel function $h(t,x)$ on $[0,
\infty)\times\bR^d$ is \textit{parabolic} with respect to the process
$X^m$ in a relatively open subset $E$ of $[0, \infty)\times\bR^d$
if for every relatively compact open subset $E_1$ of $E$,
$h(t, x)=\E^{(t,x)}
[ h ( Z^m_{\widetilde\tau^{m}_{E_1}} ) ]$
for every $(t, x)\in E_1, $ where $\widetilde\tau^m_{E_1} =\inf\{s>
0\dvtx Z^m_s\notin E_1\}$.
Note that $p_D^m(\cdot, \cdot, y)$ is parabolic with respect to the
process $X^m$.

\begin{theo}\label{T:2.3}
For any $R>0$ and $M>0$, there exists
$C_{15}>0$ such that for every
$m\in(0, M]$, $\delta\in(0, 1) $, $x_0\in\bR^d$, $t_0\ge
0$, $r\in(0,R]$ and every nonnegative function $u$ on $[0,
\infty)\times\bR^d$ that is
parabolic with respect to the
process $X^m$ on
$(t_0,t_0+4\delta r^\alpha] \times B(x_0,4r)$,
\[
\sup_{(t_1,y_1)\in Q_-}u(t_1,y_1)\le
C_{15} \inf_{(t_2,y_2)\in
Q_+}u(t_2,y_2),
\]
where $Q_-=[t_0+\delta r^\alpha,t_0+2\delta r^\alpha]\times B(x_0,r)$
and $Q_+=[t_0+3\delta r^\alpha,t_0+ 4\delta r^\alpha]\times B(x_0,r)$.
\end{theo}

\begin{pf}
Since $\psi$ is decreasing,
by the change of variable $z=|y|w$, we have for any $|y|\ge2r$,
\begin{eqnarray*}
\frac{1}{r^d}
\int_{B(0,r)} \frac{\psi(m^{1/\alpha}|z-y|)\,dz}{|z-y|^{d+\alpha}}
&\ge&
\frac{\psi(m^{1/\alpha}|y|)} {r^{d}|y|^\alpha}
\int_{ \{ |w| \le{r}/{|y|}, |w-{y}/{|y|}| \le1 \}}
\frac{dw}{|w-\fraca{y}{|y|}|^{d}}\\
&\ge&
c_0\frac{\psi(m^{1/\alpha}|y|)} {|y|^{d+\alpha}}.
\end{eqnarray*}
Thus there is a constant $c>0$ so that for every $m>0$,
\[
J^m(x, y) \leq\frac{c}{r^d} \int_{B(x, r)} J^m(z, y)\,dz
\qquad\mbox{for every } r\leq\frac{|x-y|}2.
\]
The above property is called UJS (see~\cite{CKK,CKK2}). Using
Theorem~\ref{T:2.2} and {UJS},
the conclusion of the theorem now follows from
\cite{CKK}, Theorem~4.5, or~\cite{CKK2}, Theorem~5.2.\vspace*{-2pt}
\end{pf}

\section{Preliminary lower bound estimates}\label{S:ub}

In this section, we give some preliminary lower bounds on $p^m_D(t, x,
y)$, which will be
used in Section~\ref{S:lb}
to derive the sharp two-sided estimates for $p^m(t, x, y)$ as well as
for $p^m_D(t, x, y)$.\vspace*{-2pt}

\begin{lemma}\label{L:4.2}
For any positive constants $M$, $T$, $b$ and $a$, there exists
$C_{16}=C_{16}(a, b, M, \alpha, T)>0$ such that for all $m\in(0, M]$,
$z \in\bR^d$ and $\lambda\in(0, T]$,
\[
\mathop{\mathop{\inf}_{y\in\bR^d }}_{|y -z| \le b \lambda
^{1/\alpha}}
\P_y \bigl(
\tau^m_{B(z, 2b \lambda^{1/\alpha} )} > a\lambda\bigr) \ge
C_{16}.\vspace*{-2pt}
\]
\end{lemma}

\begin{pf} By Theorem~\ref{T:2.2}, there exists $\eps=\eps(b, \alpha
, M,
T)>0$ such that for all $m\in(0, M]$ and $\lambda\in(0, T]$,
\[
\inf_{y\in\bR^d } \P_y \bigl( \tau^m_{B(y, b \lambda^{1/\alpha
}/2 )} >
\eps\lambda\bigr) \ge\frac12.
\]
We may assume that $\eps<a$. Applying Theorem~\ref{T:2.3} at most
$1+ [ (a-\eps) (4/b)^\alpha]$ times, we get that there exists
$c_1=c_1(\alpha, M, a, T)>0$ such that for all $m\in(0, M]$,
\[
c_1 p ^m_{B(y, b \lambda^{1/\alpha} )}(\eps\lambda,y,w) \le
p ^m_{B(y, b \lambda^{1/\alpha} )}(a \lambda,y,w)\qquad
\mbox{for } w \in B(y, b \lambda^{1/\alpha}/2 ).
\]
Thus for any $m\in(0, M]$,
\begin{eqnarray*}
\P_y \bigl( \tau^m_{B(y, b \lambda^{1/\alpha} )} > a\lambda\bigr)
&=& \int_{B(y, b \lambda^{1/\alpha} )}
p ^m_{B(y, b \lambda^{1/\alpha} )}(a \lambda,y,w)\,dw\\
&\ge& \int_{B(y, b \lambda^{1/\alpha}/2 )} p ^m_{B(y, b
\lambda^{1/\alpha} )}(a \lambda,y,w)\,dw \\
&\ge& c_1 \int_{B(y, b \lambda^{1/\alpha}/2 )} p ^m_{B(y, b
\lambda^{1/\alpha}/2 )}(\eps\lambda,y,w)\,dw \ge c_1/2.
\end{eqnarray*}
This proves the lemma.\vspace*{-2pt}
\end{pf}

For the next four results, $D$ is an arbitrary nonempty open set and
we use the convention that $\delta_{D}(\cdot) \equiv\infty$ when $D
=\bR^d$.\vspace*{-2pt}

\begin{prop}\label{step1}
Let $M$ and $T$ be positive constants. Suppose that $(t, x,  y)\in
(0, T]\times D\times D$ with $\delta_D(x) \ge t^{1/\alpha} \geq
2|x-y|$. Then there exists a~positive constant $C_{17}=C_{17}(M,
\alpha, T)$ such that for any $m\in(0, M]$,
%
%
\begin{equation}\label{e:lb1}
p^m_D(t,x,y) \ge C_{17} t^{-d/\alpha}.\vspace*{-2pt}
\end{equation}
\end{prop}

\begin{pf} Let $t \le T$ and $x, y \in D$ with $\delta_D(x) \ge
t^{1/\alpha} \geq
2|x-y|$.
Note that, since $t \le T$, we have
$
|x-y|
\le2^{-1}t^{1/\alpha}\le2^{-1}T^{1/\alpha}.
$
Thus,\vadjust{\goodbreak} by the uniform parabolic Harnack inequality (Theorem~\ref
{T:2.3}), there exists \mbox{$c_1=c_1(\alpha, M, T)>0$} such that for
any $m\in(0, M]$,
\[
p^m_D(t/2, x, w) \le c_1 p^m_D(t,x,y) \qquad\mbox{for every
} w \in B(x, 2t^{1/\alpha}/3) .
\]
This together with Lemma~\ref{L:4.2} yields that for any $m\in(0,
M]$,
\begin{eqnarray*}
p^m_D(t, x, y) &\geq& \frac{1}{c_1 | B(x,t^{1/\alpha}/2)|}
\int_{B(x,t^{1/\alpha}/2)} p^m_D(t/2, x, w)\,dw\\
&\geq& c_2 t^{-d/\alpha} \int_{B(x,t^{1/\alpha}/2)}
p^m_{B(x,t^{1/\alpha}/2)} (t/2, x, w)\,dw \\
&=& c_2 t^{-d/\alpha} \P_x \bigl( \tau^m_{B(x,t^{1/\alpha}/2)} >
t/2 \bigr) \geq c_3 t^{-d/\alpha},
\end{eqnarray*}
where $c_i=c_i(T, \alpha, M)>0$ for $i=2, 3$.
\end{pf}

\begin{lemma}\label{l:st3_3}
Let $M>0$ and $T>0$ be constants.
Suppose that $(t, x, y)\in(0, T]\times D\times D$ with $\min\{
\delta_D(x), \delta_D (y) \}
\ge t^{1/\alpha}$ and
$t^{1/\alpha} \leq2|x-y|$.
Then there exists a constant
$C_{18}=C_{18}(\alpha, T, M)>0$ such that for all $m\in(0, M]$,
\[
\P_x \bigl( X^{m,D}_t \in B ( y, 2^{-1} t^{1/\alpha} )
\bigr) \ge C_{18} t^{d/\alpha+1} J^m(x,y).
\]
\end{lemma}

\begin{pf}
By Lemma~\ref{L:4.2}, starting at $z\in B(y, 4^{-1}
t^{1/\alpha})$, with probability at least $c_1=c_1(\alpha, M, T)>0$,
for any $m\in(0, M]$,
the process $X^m$ does not move more than $6^{-1} t^{1/\alpha} $ by
time $t $. Thus, it is sufficient to show that there exists a
constant $c_2=c_2(\alpha, M, T)>0$ such that for any $m\in(0, M]$,
$t \in(0, T] $ and $(x, y)$ with
$t^{1/\alpha} \leq2|x-y| $,
%
%
\begin{equation}\label{eq:molow}
\P_x (X^{m,D} \mbox{ hits the ball } B(y,
4^{-1} t^{1/\alpha})\mbox{ by time } t ) \ge c_2
{t^{d/\alpha+1}}{ J^m(x,y)}.\hspace*{-35pt}
\end{equation}

Let $B_x:=B(x, 6^{-1} t^{1/\alpha})$,
$B_y:=B(y, 6^{-1} t^{1/\alpha})$
and $\tau^m_x:=\tau^m_{B_x}$. It follows from Lemma~\ref{L:4.2},
there exists $c_3=c_3(\alpha, M, T)>0$ such that for all $m\in(0,
M]$,
%
%
\begin{equation}\label{eq:lowtau}
\E_x [ t \wedge\tau^m_{x} ] \ge t \P_x
(\tau^m_{x} \ge t )
\ge c_3 t
\qquad\mbox{for }
t>0.
\end{equation}
By the L\'{e}vy system in \eqref{e:levy},
%
%
\begin{eqnarray} \label{e:dew}
&&\P_x (X^{m,D} \mbox{ hits the ball } B(y,
4^{-1} t^{1/\alpha} )\mbox{ by time } t ) \nonumber\\
&& \qquad\ge\P_x\bigl(X^m_{t\wedge\tau^m_x}\in B(y, 4^{-1}a
t^{1/\alpha})
\mbox{ and }
t \wedge\tau^m_x \mbox{ is a jumping time}\bigr)\\
&& \qquad\geq\E_x \biggl[\int_0^{t\wedge\tau^m_x} \int_{B_y}
J^m(X^m_s, u) \,du\,ds \biggr].\nonumber
\end{eqnarray}

We consider two cases separately.

(i) Suppose $|x-y| \le T^{1/\alpha}$.
Since $|x-y| \ge2^{-1} t^{1/\alpha}$, we have for
$s < \tau^m_{x}$ and $u \in B_y$,
\[
|X^m_s-u|\le|X^m_s-x|+|x-y|+|y-u| \le2|x-y|.\vadjust{\goodbreak}
\]
Thus from \eqref{e:dew}, for any $m\in(0, M]$,
\begin{eqnarray*}
&&\P_x (X^{m,D} \mbox{ hits the ball } B(y,
4^{-1} t^{1/\alpha} )\mbox{ by time } t )\\[-2pt]
&& \qquad\ge\E_x [ t\wedge\tau^m_x ] \int_{B_y}
j^m(2|x-y|)\, du\\[-2pt]
&& \qquad\ge c_4 t |B_y| j^m(2|x-y|) \ge c_5
{t^{d/\alpha+1}}{ j^m(2|x-y|)}
\end{eqnarray*}
for some positive constants $c_i=c_i(\alpha, M, T)$, $i=4, 5$. Here
in the second inequality above, we used (\ref{eq:lowtau}). Therefore
in view of \eqref{H:1},
the assertion of the lemma holds
when $|x-y| \le T^{1/\alpha}$.

(ii) Suppose $|x-y| > T^{1/\alpha}$.
In this case,
for
$s < \tau^m_{x}$ and $u \in B_y$,
\begin{eqnarray*}
|X^m_s-u|& \le& |X^m_s-x|+|x-y|+|y-u|
\\[-2pt]
&\le& |x-y| +3^{-1}
t^{1/\alpha} \le|x-y|
+3^{-1}
T^{1/\alpha}.
\end{eqnarray*}
Thus from \eqref{e:dew},
for any $m\in(0, M]$,
\begin{eqnarray*}
&&\P_x (X^{m,D} \mbox{ hits the ball } B(y,
4^{-1} t^{1/\alpha} )\mbox{ by time } t )\\[-2pt]
&& \qquad\ge\E_x [ t\wedge\tau^m_x ] \int_{B_y}
j^m (|x-y| +3^{-1}
T^{1/\alpha} )\, du\\[-2pt]
&& \qquad\ge c_6 t |B_y| j^m (|x-y| +3^{-1}
T^{1/\alpha} ) \\[-2pt]
&& \qquad\ge c_7 {t^{d/\alpha+1}}{ j^m (|x-y| +3^{-1}
T^{1/\alpha} )}
\end{eqnarray*}
for some positive constants $c_i=c_i(\alpha, M, T)$, $i=6, 7$. Here
in the second inequality, (\ref{eq:lowtau}) is used. Since $|x-y| >
T^{1/\alpha}
$, by
\eqref{H:2}, we see that the assertion of the lemma is valid for
$|x-y| > T^{1/\alpha}$ as well.
\end{pf}

\begin{prop}\label{step3}
Let $M$ and $T$ be positive constants. Suppose that $(t, x,  y)\in
(0, T]\times D\times D$ with $\min\{ \delta_D(x), \delta_D
(y) \} \ge(t/2)^{1/\alpha}$ and
$(t/2)^{1/\alpha} \leq
2|x-y|$.
Then there exists a constant $C_{19}=C_{19}(\alpha, M, T)>0$ such
that for all $m\in(0, M]$,
\[
p^m_D(t, x, y) \ge C_{19} { t}{J^m(x,y)}.
\]
\end{prop}

\begin{pf}
By the semigroup property, Proposition~\ref{step1} and Lemma
\ref{l:st3_3}, there exist positive constants
$c_1=c_1(\alpha, T, M)$ and $c_2=c_2(\alpha, T, M)$ such that for
all $m\in(0, M]$,
\begin{eqnarray*}
p^m_D(t, x, y) &=& \int_{D} p^m_{D}(t/2, x, z)
p^m_{D}(t/2, z, y)\,dz\\[-2pt]
&\ge& \int_{B(y, 2^{-1}
(t/2)^{1/\alpha})}
p^m_{D}(t/2, x,
z) p^m_{D}(t/2, z, y)\,dz\\[-2pt]
&\ge& c_1 t^{-d/\alpha} \bP_x \bigl( X^{m,D}_{t/2} \in B(y, 2^{-1}
(t/2)^{1/\alpha}) \bigr)\\[-2pt]
&\ge& c_2 { t}{J^m(x,y)}.
\end{eqnarray*}
\upqed
\end{pf}\eject

Combining Propositions~\ref{step1} and~\ref{step3}, we have the
following preliminary lower bound for $p_D^m(t,x,y)$.\vspace*{-2pt}

\begin{prop}\label{step31}
Let $M$ and $T$ be positive constants. Suppose that $(t, x,  y)\in
(0, T]\times D\times D$ with $\min\{ \delta_D(x), \delta_D
(y) \} \ge t^{1/\alpha}$.
Then there exists a~constant $C_{20}=C_{20}(\alpha, M, T)>0$ such
that for all $m\in(0, M]$,
\[
p^m_D(t, x, y) \ge C_{20} \bigl(t^{-d/\alpha} \wedge{
t}{J^m(x,y)}\bigr).
\]
\end{prop}

\section{Sharp two-sided Dirichlet heat kernel estimates} \label{S:lb}

The goal of this section is to establish the
sharp two-sided estimates for $p^m_D(t, x, y)$ as stated in
Theorem~\ref{t:main}.

First,
combining \eqref{e:ma2} and \eqref{e:mb2} with
Proposition~\ref{step31},
we have the following sharp two-sided estimates for
$p^m(t,x,y)$.

\begin{theo}\label{T:lbRd}
Let $M$ and $T$ be positive constants.
Then there exists a~constant
$C_{21}=C_{21}(\alpha, M, T)> 1$
such that for all $m\in(0, M]$, $t\in(0, T]$ and $x, y\in\bR^d$,
\[
C_{21}^{-1} \bigl( t^{-d/\alpha} \wedge{ t}{J^m(x,y)} \bigr)
\leq p^m(t, x, y) \leq C_{21} \bigl( t^{-d/\alpha} \wedge{
t}{J^m(x,y)} \bigr) .
\]
\end{theo}

The two-sided estimates in Theorem~\ref{T:lbRd} will be used
in the proof of Theorem~\ref{t:ub} to derive
sharp uniform upper bound on the Dirichlet heat kernel
$p^m_D(t, x, y)$.

\begin{lemma}\label{L:2.1}
Suppose $M>0$ and $r_0 \le R_0 M^{-1/\alpha}$. Let $E=\{x\in\bR
^d\dvtx\break
|x| >r_0\}$.
For every $T>0$, there is a
constant $C_{22}=C_{22}(r_0, \alpha, M, T)>0$ such that
\[
p^m_E(t, x, y) \leq C_{22} \sqrt{t} {\delta_E(x)^{\alpha/2}}
j^{m}( |x-y|/16)
\]
for all $m\in(0, M]$, $r_0<|x|<5r_0/4$, $ |y|\geq2r_0$ and $t\leq
T$.
\end{lemma}

\begin{pf}
Define
\[
U:=
\cases{\displaystyle\{z\in\bR^d\dvtx r_0<|z| < 3r_0/2 \},&\quad
if $d \ge2$,\cr\displaystyle
\{z\in\bR^1\dvtx r_0<z < 3r_0/2 \},&\quad if $d =1$.
}
\]
It is well known (see, e.g.,~\cite{Sz}) that $X^m_{\tau_U}\notin
\partial U$.
For $ r_0<|x|<5r_0/4$, $ |y|\geq2r_0$ and $t\in(0, T]$, it
follows from the strong Markov property and \eqref{e:levy} that
\begin{eqnarray*}
&&p^m_E(t, x, y)\\
&& \qquad= \E_x [ p^m_E(t-\tau^m_U, X^m_{\tau^m_U}, y);
\tau^m_U < t , (3r_0/4)+(|y|/2)\geq|X^m_{\tau^m_U}|>3r_0/2
]\\
&& \qquad \quad{}+ \E_x [ p^m_E(t-\tau^m_U, X^m_{\tau^m_U}, y);
\tau^m_U < t, |X^m_{\tau^m_U}|>
(3r_0/4)+(|y|/2) ]\\
&& \qquad\le \biggl(
\mathop{\mathop{\sup}_{{s\dvtx s\in(0, t)}}}_
{w\dvtx (3r_0/4)+(|y|/2)\geq|w|>3r_0/2}
p^m_E(t-s, w, y) \biggr)\\
&& \quad\qquad{} \times\P_x \bigl( \tau^m_U < t ,
(3r_0/4)+(|y|/2)\geq
|X^m_{\tau^m_U}|>3r_0/2 \bigr)
\\
&& \quad\qquad{} + \int_0^t \int_U p_U(s, x, z)\\
&& \qquad \quad\hphantom{{}+\int_0^t \int_U}{}\times\biggl( \int
_{\{ w\dvtx |w|>
(3r_0/4)+(|y|/2)\}} J^m(z, w) p^m_E(t-s, w, y)\,dw \biggr) \,dz
\,ds \\
&& \qquad=: I+\mathit{II}.
\end{eqnarray*}

If $|w|\leq(3r_0/4)+(|y|/2)$, then $|w-y| \geq|y|-|w| \geq
\frac12 (|y|-\frac{3r_0}2 )\geq\frac{|y|}8\geq
\frac{|x-y|}{16}$. Thus
by \eqref{e:ma2} and the fact that $p^m_E\le p^m$,
we have for
$|w|\leq(3r_0/4)+(|y|/2)$ and $0<s<t<T$,
\[
p_E^{m}(t-s, w, y)
\le p^{m}(t-s, x/16, y/16)
\le L {t e^{mT}}
j^{m}( |x-y|/16).
\]
Therefore
\[
I\leq L t e^{MT} j^{m}( |x-y|/16)
\P_x ( |X^m_{\tau^m_U}| > 3r_0/2 ).
\]
By Corollary~\ref{c:ufgfnest1},
\[
\P_x ( |X^m_{\tau^m_U}|> 3r_0/2 ) \le
C_{13}\P_x ( |X_{\tau_U}|> 3r_0/2 )
\le c_1\delta_U(x)^{\alpha/2}=c_1\delta_E(x)^{\alpha/2}
\]
for some positive constant $c_1=c_1(M, r_0, \alpha)$. Here the last
inequality is due to the boundary Harnack inequality for $X$ on $U$
proved in~\cite{B} (see the proof of~\cite{CKS}, Lemma 2.2). Thus we
have
%
%
\begin{equation}\label{e:2.3}
I\leq c_2
t e^{MT} { \delta_E (x)^{\alpha/2}} j^{m}({|x-y|}/{16}
), \qquad
m\in(0, M],
\end{equation}
for some positive constant $c_2=c_2(r_0, \alpha, M)$.

On the other hand, for $z\in U$ and $w\in\bR^d$ with
$|w|>(3r_0/4)+(|y|/2)$, we have
\[
|z-w| \geq|w|-|z| \geq
\frac12 \biggl(|y|-\frac{3r_0}2 \biggr) \geq\frac{|y|}8\geq
\frac{|x-y|}{16}.
\]
Thus by the symmetry of $p^m_E(t-s, w, y)$ in
$(w, y)$, we have that there exists $c_3=c_3(M, r_0, \alpha)>0$ such
that for any $m\in(0, M]$,
\begin{eqnarray*}
\mathit{II} &\leq& \int_0^t \biggl( \int_U p^m_U(s, x, z) \\
&&{} \hphantom{\int_0^t\biggl(} \times\biggl( \int_{\{ w\dvtx
|w|> (3r_0/4)+(|y|/2)\}} J^{m}(x/16, y/16)
p^m_E(t-s, y, w)\,dw \biggr)\,dz \biggr)\,ds \\
&\leq& {c_3} j^{m}(|x-y|/16) \int_0^t \biggl( \int_U p^m_U(s,
x, z)\,dz \biggr)\,ds.
\end{eqnarray*}

By \eqref{e:Meyer}, there exists $c_4=c_4(\alpha, T)>0$ such that for
every $s \le T$,
\[
p^{m}_U (s, x, z) \leq e^{ms} p_{U}(s,x,z) \leq c_{4}
e^{ms} \frac{\delta_U(x)^{\alpha/2}}{\sqrt{s}} \biggl(
s^{-d/\alpha}
\wedge\frac{s}{|x-z|^{d+\alpha}} \biggr).
\]
The last inequality above comes from~\cite{CKS}, Theorem~1.1.
Thus
\begin{eqnarray*}
&&\int_0^t \biggl( \int_U p^m_U(s, x, z)\,dz \biggr)\,ds \\
&& \qquad\le c_{4} e^{mT} \delta_U(x)^{\alpha/2}\\
&& \qquad  \quad {}\times \biggl( \int_0^t
\int_{\{|z|
\le s^{1/\alpha}\}} s^{-d/\alpha-1/2} \,dz \,ds + \int_0^t \int_{\{|z|
> s^{1/\alpha}\}} \frac{\sqrt{s}}{|z|^{d+\alpha}} \,dz \,ds \biggr
)\\
&& \qquad\le c_{5}\delta_E(x)^{\alpha/2} \sqrt{t}.
\end{eqnarray*}
This together with our estimate on $I$ above completes the proof the
lemma.
\end{pf}

Recall that an open set $D$
is said to satisfy
the weak uniform exterior ball
condition with radius $r_0>0$ if, for every
$z\in\partial D$, there is a ball $B^z$ of radius $r_0$
such that $B^z\subset
\bR^d\setminus\overline D$ and
$z \in\partial B^z$.

\begin{lemma}\label{ub11}
Let $M>0$ be a constant and $D$ an open set satisfying the
weak uniform exterior ball condition with radius
$r_0>0$.
For every
$T>0$, there exists a positive constant $C_{23}=C_{23}(T, r_0,
\alpha, M)$ such that for any $m\in(0, M]$ and $(t, x, y)\in(0,
T]\times D\times D$,
\[
p^m_D(t, x, y) \leq C_{23}
\biggl( 1 \wedge\frac{\delta_D(x)^{\alpha/2}}{\sqrt{t}} \biggr)
p^m(t, x/16, y/16).
\]
\end{lemma}

\begin{pf} Let $r_1=r_0 \wedge(R_0M^{-1/\alpha})$.
In view of Theorem~\ref{T:lbRd},
it suffices to prove the theorem for
$x\in D$ with $\delta_D (x)<r_1/4$. By \eqref{e:Meyer} and
\cite{CKS}, Theorem~1.1, there exists $c_1=c_1(\alpha, T, D)>0$
such that on $(0, T] \times D \times D$
%
%
\begin{equation}\label{e:6}
\qquad p^{m}_D (t, x, y) \leq e^{mt}
p_{D}(t,x,y)
\leq c_1 e^{MT}
\frac{\delta_D(x)^{\alpha/2}}{\sqrt{t}} \biggl( t^{-d/\alpha}
\wedge\frac{t}{|x-y|^{d+\alpha}} \biggr) .
\end{equation}
For $x, y\in D$, let $z\in\partial D$ so that $|x-z|=\delta_{D}
(x)$. Let $B_z\subset D^c$ be the ball with radius $r_1$ so that
$\partial B_z \cap\partial D=\{z\}$. When $\delta_{ D} (x) < r_1/4$
and $|x-y|\geq5r_1$, we have $\delta_{B_z^c} (y)> 2r_1 $ and so by
Lemma~\ref{L:2.1}, there is a constant $c_2=c_2(r_1, T, M,
\alpha)>0$ such that for any $m\in(0, M]$ and $(t, x, y)\in(0,
T]\times D\times D$,
%
%
\begin{eqnarray}\label{e:2.5n}
p^{m }_{ D} (t, x, y) &\leq& p^{m }_{(\overline B_z)^c}(t,x, y)
\leq c_2 \delta_{(\overline
B_z)^c}(x)^{\alpha/2} \sqrt{t} j^{m}(|x-y|/16) \nonumber
\\[-8pt]
\\[-8pt]
&=& c_2
\delta_{ D}(x)^{\alpha/2} \sqrt{t} j^{m}( |x-y|/16).
\nonumber
\end{eqnarray}
Since there exist constants $c_3$ and $c_4$ depending only on $M,
\alpha$ and $r_1$ such that
%
\begin{eqnarray}
\frac{c_3}{|x-y|^{d+\alpha}}\le j^{m}(|x-y|/16) \le
\frac{c_4}{|x-y|^{d+\alpha}} \nonumber\\
\eqntext{\mbox{for } m\in(0, M]
\mbox{ and } |x-y|<5r_1,}
\end{eqnarray}
combining \eqref{e:6} and \eqref{e:2.5n}
with Theorem~\ref{T:lbRd},
we arrive at the conclusion of
the theorem.
\end{pf}

\begin{theo}\label{t:ub}
Let $M$ and $T$ be positive constants. Suppose that $D$ is an open
set satisfying the weak uniform exterior ball condition with radius
$r_0>0$.
Then there exists a constant
$C_{24}=C_{24}(T, r_0, M, \alpha)>0$ such that for
all $m\in(0, M]$, $t\in(0, T]$ and $x, y\in D$,
%
%
\begin{equation}\label{e:1}
 p^m_D(t, x, y) \leq C_{24} \biggl( 1\wedge
\frac{\delta_D(x)^{\alpha/2}}{\sqrt{t}} \biggr)\biggl ( 1\wedge
\frac{\delta_D(y)^{\alpha/2}}{\sqrt{t}} \biggr)
p^m (t, x/16, y/16).\hspace*{-35pt}
\end{equation}
\end{theo}

\begin{pf} Fix $T>0$ and $M>0$.
By Lemma~\ref{ub11}, symmetry and the semigroup property,
we have for any $m\in(0, M]$ and $(t, x, y) \in(0, T] \times D \times D$,
\begin{eqnarray*}
p^m_D(t, x, y)
&=& \int_D p^m_D(t/2, x, z) p^m_D (t/2, z, y)\,dz\\
&\leq& c_1 \biggl(1\wedge\frac{\delta_{D}( x)^{\alpha/2}}{\sqrt{t}}
\biggr)\biggl (1\wedge\frac{\delta_{D}( y)^{\alpha/2}}{\sqrt{t}}
\biggr) \\
&& {}\times
\int_{\bR^d} p^m(t/2, x/16, z/16)p^m(t/2, z/16, y/16)\,dz \\
&\leq& c_2 \biggl(1\wedge\frac{\delta_{D}(
x)^{\alpha/2}}{\sqrt{t}} \biggr)\biggl (1\wedge\frac{\delta_{D}(
y)^{\alpha/2}}{\sqrt{t}} \biggr) p^m(t, x/16, y/16) .
\end{eqnarray*}
\upqed
\end{pf}

In the next two results,
the open set $D$ is assumed to satisfy the
uniform interior ball condition with radius $r_0>0$ in the following
sense: For every \mbox{$x\in D$} with $\delta_D (x)<r_0$, there is $z_x\in
\partial D$ so that $|x-z_x|=\delta_D(x)$ and $B(x_0, r_0)\subset D$ for
$x_0:= z_x+r_0 (x-z_x)/|x-z_x|$.
Note that this condition is strictly stronger than the
weak uniform interior ball condition with radius~$r_0$
defined as follows: For every $z\in\partial D$,
there is a ball $B^z$ of radius $r_0$ such that $B^z\subset
D$ and $z \in\partial B^z$.
Here is an example. In $\bR^2$, let $x_k=(2k, 0)\in\bR^2$ and
define $D=\bR^2\setminus\bigcup_{k=1}^\infty\partial B(x_k, 1/k)$.
Then $D$ satisfies the weak uniform interior ball condition but not
the uniform interior ball condition.\looseness=1

Under the uniform interior ball condition, we will
prove the following lower
bound for $p^m_D(t, x, y)$.

\begin{theo}\label{t:lb}
For any $M>0$ and $T>0$ there exists positive constant $C_{25} =
C_{25}(\alpha, T, M, r_0)$ such that for all $m\in(0, M]$, $(t, x,
y)\in(0, T]\times
D\times D$,
\[
p^m_D(t, x, y) \ge C_{25} \biggl( 1\wedge
\frac{\delta_D(x)^{\alpha/2}}{\sqrt{t}} \biggr)\biggl ( 1\wedge
\frac{\delta_D(y)^{\alpha/2}}{\sqrt{t}} \biggr)
\bigl(t^{-d/\alpha}\wedge{t}{j^{m}(|x-y|)} \bigr).
\]
\end{theo}

In order to prove the theorem,
for $M>0$, we let
%
%
\begin{equation}\label{tsmall}
T_0=T_0(r_0, R_0, M): = \biggl(\frac{r_0 \wedge R_0 M^{-1/\alpha
}}{16} \biggr)^{\alpha}.
\end{equation}

In the remainder of this section, for any $x \in D$ with $\delta_D(x)
< r_0$,
$z_x$ is a~point on $\partial D$ such that $|z_x-x|=\delta_D (x)$
and $\mathbf{n}(z_x):=(x-z_x)/|z_x-x|$.

\begin{lemma}\label{l:st2_3}
Let $M>0$ be a constant. Suppose that $(t, x)\in(0, T_0]\times D$
with $\delta_D(x) \leq3 t^{1/\alpha} < r_0/4$ and $\kappa\in(0, 1)$.
Put $x_0=z_x+ 4.5t^{1/\alpha}\mathbf{n}(z_x)$.
Then for
any $a>0$, there exists a constant $C_{26}=C_{26}(M, \kappa, \alpha,
r_0, a)>0$ such that for all $m\in(0, M]$,
%
%
\begin{equation}\label{e:4.5}
\P_x \bigl( X^{m, D}_{at} \in B(x_0, \kappa t^{1/\alpha}) \bigr)
\ge C_{26} \frac{ \delta_D(x)^{\alpha/2}}{t^{1/2}}.
\end{equation}
\end{lemma}

\begin{pf} Let $0< \kappa_1 \le\kappa$ and assume first that
$2^{-4}\kappa_1 t^{1/\alpha} < \delta_D(x)\leq3t^{1/\alpha}$. As in
the proof of Lemma~\ref{l:st3_3}, we get that, in this case, using
the fact that $|x-x_0|\in[1.5\kappa t^{1/\alpha}, 6t^{1/\alpha}]$,
there exist constants $c_i=c_i(\alpha, \kappa_1, M, r_0, a)>0$,
$i=1,2,$ such that for all $m\in(0, M]$ and $t\le T_0$,
%
%
\begin{equation}\label{e:case1}
\P_x \bigl( X^{m, D}_{at} \in B(x_0, \kappa_1 t^{1/\alpha}) \bigr)
\ge c_1t^{d/\alpha+1}J^m(x, x_0) \ge c_2>0.
\end{equation}
By taking
$\kappa_1=\kappa$, this shows that \eqref{e:4.5} holds for all $a>0$
in the case when
$2^{-4}\kappa t^{1/\alpha} < \delta_D(x)\leq3t^{1/\alpha}$.

So it suffices to consider the case that $\delta_D(x) \leq
2^{-4}\kappa t^{1/\alpha}$. We now show that there is some $a_0>1$
so that \eqref{e:4.5} holds for every $a\geq a_0$ and $\delta_D(x)
\leq2^{-4}\kappa t^{1/\alpha}$. For simplicity, we assume without
loss of generality that $x_0=0$ and let $\wh B:=B(0, \kappa
t^{1/\alpha})$. Let
$x_1
=z_x+4^{-1}\kappa\mathbf{n}(z_x)t^{1/\alpha}$
and $B_1:= B(x_1, 4^{-1}\kappa t^{1/\alpha})$.
By the strong Markov property of $X^{m, D}$ at the first exit time
$\tau^m_{B_1}$ from~$B_1$ and Lemma~\ref{L:4.2}, there exists
$c_3=c_3(a, \kappa, \alpha, M, T)>0$ such that for all $m\in(0,
M]$,
%
%
\begin{eqnarray}
\label{e:4.8}
&&\P_{x} (X^{m, D}_{at} \in{\wh B} )\nonumber\\
&& \qquad\geq\P_{x} \bigl( \tau^m_{B_1}<at, X^m_{\tau^m_{B_1}}
\in B(0,
2^{-1}\kappa t^{1/\alpha})
\mbox{ and}\nonumber
\\[-8pt]
\\[-8pt]
&&\hphantom{ \qquad \quad\P_{x} \bigl(}|X^{m, D}_s-X^m_{\tau
_{B_1}}|<2^{-1}\kappa t^{1/\alpha} \mbox{ for }
s\in[\tau^m_{B_1}, \tau^m_{B_1}+at] \bigr) \nonumber\\
&& \qquad\geq c_3 \P_{x} \bigl( \tau^m_{B_1}<at \mbox{ and }
X^m_{\tau^m_{B_1}}\in B(0, 2^{-1}\kappa t^{1/\alpha}) \bigr).
\nonumber
\end{eqnarray}

It follows from the first display in Theorem~\ref{t:ufgfnest}
and the explicit formula
for the Poisson kernel of balls with respect to $X$ that there exist
\mbox{$c_4=c_4(\alpha, M)>0$} and $c_5=c_5(\alpha, M, \kappa, r_0)>0$ such
that for all $m\in(0, M]$,
%
%
\begin{eqnarray}\label{e:case4}
\P_{x} \bigl(X^m_{\tau^m_{B_1}}\in B(0, 2^{-1}\kappa
t^{1/\alpha}) \bigr)
&\ge& c_4
\P_{x} \bigl(X_{\tau_{B_1}}\in B(0,
2^{-1}\kappa t^{1/\alpha}) \bigr)\nonumber
\\[-8pt]
\\[-8pt]
&\ge& c_5
\biggl(\frac{\delta_{D}(x)}{t^{1/\alpha}} \biggr)^{\alpha
/2}.
\nonumber
\end{eqnarray}
Applying Theorem~\ref{t:ufgfnest} and the estimates for
$G_{B_1}$ (see,
e.g.,~\cite{CS1}, (1.4)), we get that there exist
$c_6=c_6(\alpha, M)>0$ and $c_7=c_7(\alpha, M, \kappa, r_0)>0$ such
that for all $m\in(0, M]$,
\begin{eqnarray*}
\P_{x}( \tau^m_{B_1}\geq a t ) \leq (at)^{-1} \E_{x} [
\tau^m_{B_1} ] \le c_6 (at)^{-1} \E_{x} [ \tau_{B_1} ] \le
a^{-1} c_7 \biggl(\frac{\delta_{D}(x)}{t^{1/\alpha}}
\biggr)^{\alpha/2}.
\end{eqnarray*}
Define $a_0= 2c_7/(c_5) $. We have by
\eqref{e:4.8} and \eqref{e:case4}
and the display above that for $a\geq a_0$ and $m\in
(0, M]$,
%
%
\begin{eqnarray}\label{e:case2}\qquad
\P_x ( X^{m, D}_{at}\in\wh B ) &\geq& c_3 \bigl ( \P_{x} \bigl(
X^m_{\tau^m_{B_1}}\in B(0, 2^{-1}\kappa t^{1/\alpha} ) \bigr) - \P_{x}
( \tau^m_{B_1}\geq at ) \bigr)\nonumber
\\[-8pt]
\\[-8pt]
&\geq& c_3 (c_5/2)
\biggl(\frac{\delta_{D}(x)}{t^{1/\alpha}}
\biggr)^{\alpha/2}.
\nonumber
\end{eqnarray}
Equations \eqref{e:case1} and \eqref{e:case2} show that \eqref{e:4.5} holds
for every $a\geq a_0$ and for every $x\in D$ with $\delta_D(x) \leq
3 t^{1/\alpha}$.

Now we deal with the case $0<a<a_0$ and $\delta_D(x) \leq
2^{-4}\kappa t^{1/\alpha}$. If $\delta_D(x) \leq3
(at/a_0)^{1/\alpha}$, we have from \eqref{e:4.5} for the case of
$a=a_0$ that there exist $c_8=c_8(\kappa, \alpha, M)>0$ and
$c_9=c_9(\kappa, \alpha, M, a)>0$ such that for all $m\in(0, M]$,
\begin{eqnarray*}
\P_x \bigl( X^{m, D}_{at}\in B(x_0, \kappa t^{1/\alpha})
\bigr) &\ge& \P_x \bigl( X^{m, D}_{a_0 (at/a_0)}\in B(x_0,
\kappa(at/a_0)^{1/\alpha}) \bigr) \\
&\geq& c_8 \biggl(\frac{\delta_D(x)}{(at/a_0)^{1/\alpha}}
\biggr)^{\alpha/2}= c_9
\biggl(\frac{\delta_{D}(x)}{t^{1/\alpha}} \biggr)^{\alpha/2}.
\end{eqnarray*}
If $3 (at/a_0)^{1/\alpha} < \delta_D(x) \leq2^{-4}\kappa
t^{1/\alpha}$ [in this case $\kappa> 3 \cdot2^4
(a/a_0)^{1/\alpha}$], we get~\eqref{e:4.5} from \eqref{e:case1} by
taking $\kappa_1=(a/a_0)^{1/\alpha}$. The proof of the lemma is now
complete.
\end{pf}

\begin{pf*}{Proof of Theorem~\ref{t:lb}}
In the first part of this proof, we
adapt
some arguments from~\cite{BGR}.

Assume first that $t \le T_0$.
Since $D$ satisfies the uniform interior ball condition
with radius $r_0$ and $0<t \le T_0$, we can choose $\xi^t_x$ as follows:\vspace*{2pt}
if $\delta_D(x) \le3 t^{1/\alpha}$,
let $\xi^t_x=z_x+ (9/2)t^{1/\alpha}\mathbf{n}(z_x)$ [so that $B(\xi
^t_x, (3/2) t^{1/\alpha})\subset B(
z_x+\break 3t^{1/\alpha}\mathbf{n}(z_x), 3t^{1/\alpha}) \setminus\{x\}$
and $\delta_D(z) \ge
3t^{1/\alpha}$ for every $z \in B(\xi^t_x, (3/2) t^{1/\alpha})$].
If $\delta_D(x) > 3 t^{1/\alpha}$, choose $\xi^t_x \in B(x, \delta_D(x))$
so that $|x-\xi^t_x|=(3/2) t^{1/\alpha}$. Note that in this case,
$B(\xi^t_x, (3/2) t^{1/\alpha})\subset B(x, \delta_D(x))
\setminus\{
x\}$ and $\delta_D(z) \ge t^{1/\alpha}$
for every $z \in B(\xi^t_x, 2^{-1} t^{1/\alpha})$.
We also define $\xi^t_y$ the same way.

If $\delta_D(x) \le3 t^{1/\alpha}$, by
Lemma~\ref{l:st2_3} (with $a=3^{-1}, \kappa=2^{-1}$),
\[
\P_x \bigl(X^{m, D}_{t/3} \in B(\xi^t_x, 2^{-1}t^{1/\alpha} )
\bigr) \ge c_0 \frac{\delta_D(x)^{\alpha/2}}{\sqrt{t}}.
\]
If $\delta_D(x) > 3 t^{1/\alpha}$, by Proposition~\ref{step31},
%
%
\begin{eqnarray}
\label{e:loww_0}
\P_x \bigl(X^{m, D}_{t/3} \in B(\xi^t_x, 2^{-1}t^{1/\alpha}
) \bigr)&=&\int_{B(\xi^t_x, 2^{-1}t^{1/\alpha} )}
p^{m}_{{D}}(t/3,x,u)\,du \nonumber\\
&\ge& c_1 t^{-d/\alpha}\bigl(1 \wedge\psi((MT)^{1/\alpha})
\bigr)
|B(\xi^t_x, 2^{-1}t^{1/\alpha} )|
\\
&\ge& c_2 \ge c_3 \biggl(\frac{\delta_D(x)^{\alpha/2}}{\sqrt{t}}
\wedge1 \biggr).\nonumber
\end{eqnarray}
Similarly,
%
%
\begin{equation}
\label{e:loww_01}
\P_y \bigl(X^{m, D}_{t/3} \in B(\xi^t_y, 2^{-1}t^{1/\alpha} )
\bigr) \ge c_3
\biggl(\frac{\delta_D(y)^{\alpha/2}}{\sqrt{t}} \wedge1 \biggr).
\end{equation}
Note that by the semigroup
property, Proposition~\ref{step31} and \eqref{e:loww_0} and \eqref
{e:loww_01},
%
%
\begin{eqnarray}\label{e:loww1}
&&p^{m}_{D}(t,x,y)\nonumber\hspace*{-35pt}\\
&& \quad\geq \int_{B(\xi^t_y, 2^{-1}t^{1/\alpha} )}\int_{B(\xi^t_x,
2^{-1}t^{1/\alpha} )}
p^{m}_{D}(t/3,x,u)
p^{m}_{D}(t/3,u,v)\nonumber\hspace*{-35pt}\\
&& \quad\hphantom{\geq \int_{B(\xi^t_y, 2^{-1}t^{1/\alpha}
)}\int_{B(\xi^t_x,
2^{-1}t^{1/\alpha} )}} {}\times p^{m}_{D}
(t/3,v,y)\,du\,dv \nonumber\hspace*{-35pt}\\[-8pt]
\\[-8pt]
&& \quad\geq c_4\int_{B(\xi^t_y, 2^{-1}t^{1/\alpha} )}
\int_{B(\xi^t_x,
2^{-1}t^{1/\alpha} )}
p^{m}_{{D}}(t/3,x,u)\bigl(tJ^{m}(u,v)\wedge
t^{-d/\alpha}\bigr)\nonumber\hspace*{-35pt}
\\
&&\quad\hphantom{\geq c_4\int_{B(\xi^t_y, 2^{-1}t^{1/\alpha} )}
\int_{B(\xi^t_x,
2^{-1}t^{1/\alpha} )}}{} \times p^{m}_{D}(1/3,v,y)\,du\,dv
\nonumber\hspace*{-35pt}\\
&& \quad\geq c_5 \biggl(\mathop{\mathop{\inf}
_{u \in B(\xi^t_x, 2^{-1}t^{1/\alpha} ) }}_{
v \in
B(\xi^t_y,2^{-1}t^{1/\alpha} )} \bigl(tJ^{m}(u,v)\wedge
t^{-d/\alpha}\bigr) \biggr) \biggl(\frac{\delta_D(x)^{\alpha/2}}
{\sqrt
{t}} \wedge1 \biggr)
\biggl(\frac{\delta_D(y)^{\alpha/2}}{\sqrt{t}} \wedge1
\biggr).\nonumber\hspace*{-35pt}
\end{eqnarray}
For $(u,v) \in B(\xi^t_x, 2^{-1}t^{1/\alpha} ) \times
B(\xi^t_y,2^{-1}t^{1/\alpha} )$, since
$|u-v| \leq t^{1/\alpha} +|\xi^t_x-\xi^t_y|\leq10t^{1/\alpha}
+|x-y|$,
by considering
the cases $|x-y| \ge t^{1/\alpha}$ and $|x-y| <t^{1/\alpha}$
separately using \eqref{H:1} and \eqref{H:2},
we have
%
%
\begin{eqnarray}\label{e:loww2}
&&\inf_{(u,v) \in B(\xi^t_x, 2^{-1}t^{1/\alpha} ) \times B(\xi
^t_y,2^{-1}t^{1/\alpha} )}
\bigl(tJ^{m}(u,v)\wedge
t^{-d/\alpha}\bigr)\nonumber
\\[-8pt]
\\[-8pt]
&& \qquad\ge c_6 \bigl(tJ^{m}(x,y)\wedge
t^{-d/\alpha}\bigr).
\nonumber
\end{eqnarray}
Thus combining \eqref{e:loww1} and \eqref{e:loww2},
we conclude that for $t \in(0, T_0]$,
%
%
\begin{equation}\label{e:w2}
 p^m_D(t, x, y)\ge c_7 \biggl(\frac{\delta_D(x)^{\alpha
/2}}{\sqrt
{t}} \wedge1 \biggr)
\biggl(\frac{\delta_D(y)^{\alpha/2}}{\sqrt{t}} \wedge1 \biggr
)\bigl(
tJ^{m}(x,y)\wedge t^{-d/\alpha}\bigr) .\hspace*{-35pt}
\end{equation}

Next assume $T=2T_0$.
Recall that $T_0=((r_0\wedge R_0M^{-1/\alpha})/16)^\alpha$.
For $(t, x,
y)\in(T_0, 2T_0] \times D \times D$, let $x_0, y_0\in D$ be such
that $\max\{|x-x_0|, |y-y_0|\} <r_0$ and $\min\{\delta_D (x_0),
\delta_D (y_0)\}\geq r_0/2$. Note that, if $|x-y| \ge4r_0$, then
$|x-y|-2r_0 \le|x_0-y_0| \le|x-y|+2r_0$, so by \eqref{H:2},
$c_8^{-1} J^m(x_0, y_0) \le J^m(x, y) \le c_8 J^m(x_0, y_0)$ for
some constant $c_8=c_8(M)>1$. Thus
we have
%
%
\begin{equation}\label{e:N1}
(t/2)^{-d/\alpha}
\wedge\frac{t}2 J^m(x_0, y_0)
\geq c_9
\bigl( t^{-d/\alpha} \wedge tJ^m(x, y) \bigr) .
\end{equation}
Similarly, there is a positive constant $c_{10}$ such that for every
$z,w \in D$,
%
%
\begin{eqnarray} \label{e:N2}\qquad
(t/3)^{-d/\alpha} \wedge(t/3) J^m(x, z)
&\ge& c_{10} \biggl( (t/12)^{-d/\alpha} \wedge
\frac{t}{12} J^m(x_0,z) \biggr), \nonumber
\\[-8pt]
\\[-8pt]
(t/3)^{-d/\alpha} \wedge
({t}/{3})J^m(w,y)
&\ge& c_{10} \biggl( (t/12)^{-d/\alpha} \wedge
\frac{t}{12}J^m(w,y_0) \biggr).
\nonumber
\end{eqnarray}
By \eqref{e:N2}
and the lower bound estimate in Theorem
\ref{t:lb} for $p^m_D$ on $(0, T_0]\times D\times D$, we have
\begin{eqnarray*}
&&\hspace*{-4pt}p^m_D(t, x, y)\\
&&\hspace*{-4pt} \qquad = \int_{D\times D}
p^m_D(t/3, x, z)p^m_D (t/3, z, w) p^m_D(t/3, w, y) \,dz\, dw \\
&&\hspace*{-5pt} \qquad\geq c_{11} \biggl(1\wedge
\frac{\delta_D(x)^{\alpha/2}}{\sqrt{t/3}} \biggr)\biggl ( 1\wedge
\frac{\delta_D(y)^{\alpha/2}}{\sqrt{t/3}} \biggr) \int_{D\times D}
\bigl( (t/3)^{-d/\alpha} \wedge(t/3) J^m(x, z) \bigr)
\\
&&\hspace*{-4pt} \qquad \quad{} \times
\biggl(1\wedge
\frac{\delta_D(z)^{\alpha/2}}{\sqrt{t/3}} \biggr) p^m_D(t/3, z, w)
\biggl( (t/3)^{-d/\alpha} \wedge\frac{t}{3}J^m(w,y) \biggr)
\\
&&\hspace*{-4pt} \qquad \quad{} \times
\biggl(1\wedge
\frac{\delta_D(w)^{\alpha/2}}{\sqrt{t/3}} \biggr) \,dz\,dw \\
&&\hspace*{-4pt} \qquad\geq c_{12} \biggl(1\wedge
\frac{\delta_D(x)^{\alpha/2}}{\sqrt{t}} \biggr)\biggl ( 1\wedge
\frac{\delta_D(y)^{\alpha/2}}{\sqrt{t}} \biggr) \int_{D\times D}
\biggl( \bigl(t/12\bigr)^{-d/\alpha} \wedge
\frac{t}{12} J^m(x_0,z) \biggr)
\\
&&\hspace*{-4pt} \qquad \quad{} \times
\biggl(1\wedge
\frac{\delta_D(z)^{\alpha/2}}{\sqrt{t/3}} \biggr) p^m_{D}(t/3, z, w)
\biggl( (t/12)^{-d/\alpha} \wedge
\frac{t}{12}J^m(w,y_0) \biggr)\\
&&\hspace*{-4pt} \qquad \quad{} \times
\biggl(1\wedge
\frac{\delta_D(w)^{\alpha/2}}{\sqrt{t/3}} \biggr)\,dz\,dw
\end{eqnarray*}
for some positive constants $c_i, i=11, 12$.
Let $D_1:=\{z \in D\dvtx \delta_D(z) > r_0/4\}$. Clearly, $x_0,
y_0\in
D_1$ and
%
%
\begin{equation}\label{e:nem}
\min\{\delta_{D_1} (x_0), \delta_{D_1} (y_0)\}\geq
r_0/4=4(T_0)^{1/\alpha} \ge4 (t/2)^{1/\alpha}.
\end{equation}
We have by Theorem~\ref{T:lbRd}, \eqref{e:N1} and Lemma~\ref{ub11}
that
\begin{eqnarray*}
&& \int_{D\times D}
\biggl( \bigl(t/(12)\bigr)^{-d/\alpha} \wedge
\frac{t}{12} J^m(x_0,z) \biggr)
\biggl(1\wedge
\frac{\delta_D(z)^{\alpha/2}}{\sqrt{t/3}} \biggr) \\
&& \qquad \quad{} \times p^m_{D}(t/3, z, w)
\biggl( (t/12)^{-d/\alpha} \wedge
\frac{t}{12}J^m(w,y_0) \biggr)
\biggl(1\wedge
\frac{\delta_D(w)^{\alpha/2}}{\sqrt{t/3}} \biggr)\,dz\,dw\\
&& \qquad\geq c_{13} \int_{D_1 \times D_1 }
p^m_{D_1} (t/12, x_0, z) p^m_{D_1}
(t/3, z, w)p^m_{D_1} \bigl(t/(12), w,
y_0\bigr) \,dz\,dw \\
&& \qquad= c_{13}
p^m_{D_1}(t/2, x_0, y_0) \geq c_{14}
\biggl( (t/2)^{-d/\alpha}
\wedge\frac{t}2 J^m(x_0, y_0) \biggr)
\\ && \qquad\ge
c_{15}
\bigl( t^{-d/\alpha} \wedge tJ^m(x, y) \bigr)
\end{eqnarray*}
for some positive constants $c_i, i=13, \ldots, 15$. Here Proposition
\ref{step31} is used in the third inequality in view of
\eqref{e:nem}.
Iterating the above argument one can deduce that Theorem~\ref{t:lb}
holds for $T=kT_0$ for any integer $k\geq2$.
This completes the proof of the theorem.
\end{pf*}

\begin{pf*}{Proof of Theorem~\ref{t:main}}
Theorem~\ref{t:main}(i) is a combination of Theorems~\ref{t:ub} and
\ref{t:lb}, so we only need to prove Theorem~\ref{t:main}(ii).

Let $D$ be a bounded $C^{1,1}$ open set in $\bR^d$ with
$C^{1,1}$ characteristics $(r_0, \Lambda_0)$.
Clearly there is a ball $B \subset D$ whose radius depends only on
$r_0$ and $\Lambda_0$.
For each $m\geq0$, the semigroup of $X^{m, D}$ is
Hilbert--Schmidt as, by Theorem~\ref{t:main}(i)
\[
\int_{D\times D}p^m_D (t, x, y)^2 \,dx \,dy= \int_Dp^m_D(2t, x, x)\,dx
\le C_1
(2t)^{-d/\alpha}
|D|
<\infty,
\]
and hence is compact.
Let $ \{\lambda^{\alpha, m, D}_k, k=1, 2, \ldots\}$ be
the eigenvalues of $ ((m^{2/\alpha}-\Delta)^{\alpha/2}
-m)|_D$, arranged in increasing order and repeated according
to multiplicity, and let
$\{\phi^{\alpha, m, D}_k, k=1, 2, \ldots\} $ be the
corresponding eigenfunctions normalized to have unit $L^2$-norm on
$D$.
It is well known that
$\lambda^{\alpha, m, D}_1$ is strictly positive and simple,
and that $\phi^{\alpha, m, D}_1 $ can be chosen to be
strictly positive on $D$, and that $\{\phi^{\alpha, m, D}_k \dvtx
k=1, 2,
\ldots\} $
forms an orthonormal basis of $L^2(D; dx)$.

We also let $ \{\lambda^{\alpha, m, B}_k
\dvtx k=1, 2 ,\ldots\}$ be the eigenvalues of
$ ((m^{2/\alpha}-\Delta)^{\alpha/2}-m)|_B$, arranged in
increasing order and repeated
according to multiplicity.
From the domain monotonicity of the first eigenvalue,
it is easy to see that
$\lambda^{\alpha, m, B}_1 \ge\lambda^{\alpha, m, D}_1$.
Thus, using~\cite{CS8}, Theorem~3.4, we have that
for every $m \in(0, M],$
%
%
\begin{equation}\label{e:4.1}
\lambda^{\alpha, m, D}_1 \leq\lambda^{\alpha, m, B}_1 \leq
(\lambda^{B}_1+m^{2/\alpha})^{\alpha/2}-m
\leq(\lambda^{B}_1+M^{2/\alpha})^{\alpha/2}
=: c_1,\hspace*{-35pt}
\end{equation}
where $\lambda^{B}_1$ the first eigenvalue of $-\Delta|_B$.
Moreover, by the Cauchy--Schwarz inequality,
%
%
\begin{equation}\label{e:4.2}
\int_D \bigl(1 \wedge\delta_D(x)^{\alpha/2} \bigr)\phi^{\alpha
, m, D}_1 (x)\,dx \le\biggl(\int_D \bigl(1 \wedge\delta
_D(x)^{\alpha
} \bigr)\,dx \biggr)^{1/2} =: c_2.\hspace*{-35pt}
\end{equation}

Since $p^m_D (t, x, y)$ admits the following eigenfunction
expansion:
\[
p^m_D (t, x, y) =\sum_{k=1}^\infty e^{-t \lambda^{\alpha, m, D}_k}
\phi^{\alpha, m,
D}_k (x) \phi^{\alpha, m, D}_k(y)
\qquad\mbox{for } t>0 \mbox{ and } x, y \in D,
\]
we have
%
%
\begin{eqnarray}\label{e:4.3}
&&\int_{D\times D} \bigl(1 \wedge\delta_D(x)^{\alpha/2}
\bigr)p^m_D(t, x, y) \bigl(1 \wedge\delta_D(y)^{\alpha/2} \bigr
)\, dx\,
dy\nonumber
\\[-8pt]
\\[-8pt]
&& \qquad=
\sum_{k=1}^\infty e^{-t \lambda^{\alpha, m, D}_k} \biggl( \int_D
\bigl(1\wedge\delta_D(x)^{\alpha/2}\bigr)
\phi^{\alpha, m, D}_k (x)\,dx \biggr)^2.
\nonumber
\end{eqnarray}
Consequently, using the fact that $\{\phi^{\alpha, m, D}_k \dvtx k=1, 2,
\ldots\} $
forms an orthonormal basis of $L^2(D; dx)$, we have
%
%
\begin{eqnarray}\label{e:4.4}
&&\int_{D\times D} \bigl(1 \wedge\delta_D(x)^{\alpha/2}
\bigr)p^m_D(t, x, y) \bigl(1 \wedge\delta_D(y)^{\alpha/2} \bigr
)\, dx \,dy
\nonumber
\\[-8pt]
\\[-8pt]
&& \qquad\leq e^{-t \lambda^{\alpha, m, D}_1} \int_D
\bigl(1 \wedge\delta_D(x)^{\alpha} \bigr)\,dx
\nonumber
\end{eqnarray}
for all $m>0$ and $t>0$.
On the other hand, since
\[
\phi^{\alpha, m, D}_1 (x) = e^{\lambda^{\alpha, m, D}_1} \int
_Dp^m_D (1, x, y)
\phi^{\alpha, m, D}_1 (y)\,dy,
\]
by the upper bound estimate in Theorem~\ref{t:main}(i) and
\eqref{e:4.2}, we see that
for every $m \in(0, M]$ and $x\in D$,
\begin{eqnarray*}
\phi^{\alpha, m, D}_1 (x)&\leq& e^{\lambda^{\alpha, m, D}_1}
C_1 \bigl(1 \wedge\delta_D(x)^{\alpha/2} \bigr) \int_D \bigl(1
\wedge\delta_D(y)^{\alpha/2} \bigr)\phi^{\alpha, m, D}_1 (y)\,dy
\\
&\leq& e^{\lambda^{\alpha, m, D}_1} c_2 C_1 \bigl(1 \wedge\delta
_D(x)^{\alpha/2} \bigr).
\end{eqnarray*}
Hence
\begin{eqnarray*}
&& \int_D \bigl(1 \wedge\delta_D(x)^{\alpha/2} \bigr)\phi
^{\alpha, m, D}_1(x)\,dx
\\&& \qquad\geq e^{-\lambda^{\alpha, m, D}_1} (c_2 C_1)^{-1}
\int_D \phi^{\alpha, m, D}_1(x)^2 \,dx=e^{-\lambda^{\alpha, m, D}_1}
(c_2 C_1)^{-1}.
\end{eqnarray*}
It now follows from \eqref{e:4.3} that for every
$m \in(0, M]$ and $t>0$
%
%
\begin{eqnarray}
\label{e:3.36}
&& \int_{D\times D} \bigl(1 \wedge\delta_D(x)^{\alpha/2}
\bigr)p^m_D(t, x, y) \bigl(1 \wedge\delta_D(y)^{\alpha/2} \bigr
)\, dx\,
dy \nonumber\\
&& \qquad\geq e^{-t \lambda^{\alpha, m, D}_1} \biggl(\int_D \bigl(1
\wedge\delta_D(x)^{\alpha/2} \bigr)\phi^{\alpha, m, D}_1(x)\,dx
\biggr)^2 \\
&& \qquad\geq e^{-(t+2)\lambda^{\alpha, m, D}_1} (c_2
C_1)^{-2}.\nonumber
\end{eqnarray}

It suffices to prove Theorem~\ref{t:main}(ii) for $T\geq
3$. For $t\geq T$ and $x, y\in D$, observe that
%
%
\begin{equation}\label{e:4.6}
\quad p^m_D (t, x, y)=\int_{D\times D}
p^m_D (1, x, z) p^m_D(t-2, z, w)
p^m_D(1, w, y) \,dz\, dw .
\end{equation}
Since $D$ is bounded, we have by the upper bound estimate in
Theorem~\ref{t:main}(i), \eqref{e:4.1} and \eqref{e:4.4} that for
every $m \in(0, M]$,
$t\geq T$ and $x, y\in D$,
\begin{eqnarray*}
&&p^m_D(t, x, y)\\
&& \qquad\leq C_1^2 \bigl(1 \wedge\delta_D(x)^{\alpha/2} \bigr
)\bigl(1
\wedge\delta_D(y)^{\alpha/2} \bigr)
\\
&& \qquad \quad{} \times
\int_{D\times D}
\bigl(1 \wedge\delta_D(z)^{\alpha/2} \bigr) p^m_D (t-2, z, w)
\bigl(1 \wedge\delta_D(w)^{\alpha/2} \bigr)\, dz\, dw \\
&& \qquad \quad{}\leq C_1^2 \bigl(1 \wedge\delta_D(x)^{\alpha/2}
\bigr)
\bigl(1 \wedge\delta_D(y)^{\alpha/2} \bigr) e^{-(t-2)\lambda
^{\alpha,
m, D}_1}
\int_D 1\wedge\delta_D(x)^\alpha\, dx \\
&& \qquad \quad\leq c_3
\delta_D(x)^{\alpha/2} \delta_D(y)^{\alpha/2} e^{t\lambda
^{\alpha, m, D}_1}.
\end{eqnarray*}
Similarly, by the lower bound estimate in Theorem
\ref{t:main}(i) and
\eqref{e:3.36}
that for every $m \in(0, M]$, $t\geq T$ and $x, y\in D$,
\begin{eqnarray*}
&& p^m_D(t, x, y) \\
&& \qquad\ge
c_4 \bigl(1 \wedge\delta_D(x)^{\alpha/2} \bigr)\bigl(1 \wedge
\delta_D(y)^{\alpha/2} \bigr)
\\
&& \qquad \quad{} \times
\int_{D\times D}
\bigl(1 \wedge\delta_D(z)^{\alpha/2} \bigr) p^m_D (t-2, z, w)
\bigl(1 \wedge\delta_D(w)^{\alpha/2} \bigr) \,dz\, dw\\
&& \qquad\ge c_5
\delta_D(x)^{\alpha/2} \delta_D(y)^{\alpha/2} e^{-t\lambda^{\alpha
, m, D}_1}.
\end{eqnarray*}
This establishes Theorem~\ref{t:main}(ii).
\end{pf*}

\begin{remark}\label{R:4.7}
(i) In this paper, we do not use the boundary Harnack
inequality for $X^m$.
The boundary decay rate is obtained by comparing the Green function
of $X^m$ in balls and annulus with that of $X^0$ through drift transform
(see Theorem~\ref{t:ufgfnest}).

(ii) Let $Y$ be a relativistic stable-like process on $\bR^d$,
as studied in~\cite{CKK3}.
If one can establish scale invariant boundary Harnack inequality
for $Y$ in bounded $C^{1,1}$ open sets with explicit boundary decay
rate $\delta_D(x)^{\alpha/2}$, then one can easily modify the
approach of this paper to show that Theorem~\ref{t:main} holds for $Y$ with
$\phi(m^{1/\alpha}|x-y|)$ and $\phi(m^{1/\alpha}|x-y|/16)$
being replaced by $\phi(c_1|x-y|)$ and $\phi(c_2|x-y|)$, respectively,
for some positive constant $c_1$ and~$c_2$.

(iii)
By integrating \eqref{e:1.2} with
respect
to $y$,
we see that for each fixed $M, T>0$,
\[
\P_x (t<\tau^m_D)
\asymp1 \wedge\frac{\delta_D(x)^{\alpha/2}}{\sqrt{t}}
\qquad\mbox{for } m\in(0, M] \mbox{ and } (t, x)\in(0, T]\times D.
\]
Hence both \eqref{e:1.2} and \eqref{e:1.3} can be restated
as, for each fixed $T>0$ and every $(t, x, y)\in(0, T]\times D\times D$,
\begin{eqnarray*}
&& \frac1{C_1} \P_x(t<\tau^m_D) \P_y(t<\tau
^m_D)
p^m(t, x, y) \\
&& \qquad\le p^m_D(t, x, y) \le C_1
\P_x(t<\tau^m_D) \P_y(t<\tau^m_D) p^m(t, x/16, y/16).
\end{eqnarray*}
It is possible to establish the above estimates by
adapting
the approach in~\cite{BGR},
using the scale invariant boundary Harnack inequality for
$X^m$, uniform on $m \in(0, M]$, in arbitrary $\kappa$-fat open sets
which is recently established in~\cite{CKS3}, Theorem~2.6.
We omit the details here.
\end{remark}

\section*{Acknowledgments}
We thank the two referees for their
helpful comments on the first version of this paper. We also thank
Qiang Zeng
for pointing out
some typos
in the first version of this paper.

%

\printaddresses

\end{document}